\documentclass[12pt,reqno,twoside]{amsart}
\usepackage{amssymb}

\providecommand{\notinjournal}[1]{#1}

\notinjournal{
 \numberwithin{equation}{section}
 \theoremstyle{definition}
 \newtheorem*{ntn}{Notation} }

\providecommand{\toas}{\xrightarrow[N\to\infty]{}}

\def\beq#1\eeq{\begin{equation}#1\end{equation}}
\newcommand{\gap}{\vspace{1ex}}
\renewcommand{\hbar}{{\mathchar'26\mkern-9muh}}
\newcommand{\A}{\mathcal{A}}
\newcommand{\C}{\mathcal{C}}
\newcommand{\M}{\mathcal{M}}
\newcommand{\co}{\mathbb{C}}
\newcommand{\cs}{C${}^*$}
\newcommand{\R}{\mathbb{R}}
\newcommand{\Rplus}{\R_{\mkern1mu+}}
\newcommand{\Z}{\mathbb{Z}}
\newcommand{\onto}{\twoheadrightarrow}
\newcommand{\into}{\hookrightarrow}
\renewcommand{\AA}{\mathbb{A}}
\newcommand{\I}{\mathcal{I}}
\newcommand{\dlim}{\varinjlim}
\newcommand{\ilim}{\varprojlim}
\newcommand{\veclim}{\text{Vec-}\!\dlim}
\newcommand{\G}{\Gamma}
\newcommand{\Gc}{\G_{\!\mathrm{c}}}
\newcommand{\Gb}{\G_{\!\mathrm{b}}}
\newcommand{\Gn}{\G_{\!0}}
\newcommand{\Gp}{\G_{\!\text{poly}}}
\newcommand{\Lie}{\mathcal{L}}
\newcommand{\g}{\mathfrak{g}}
\newcommand{\h}{\mathfrak{h}}
\newcommand{\inner}{\mathbin{\lrcorner}}
\DeclareMathOperator{\Ad}{Ad}
\newcommand{\Orl}{\mathcal{O}_\Lambda}
\DeclareMathOperator{\End}{End}
\newcommand{\norm}[1]{\left\Vert #1 \right\Vert}
\DeclareMathOperator{\tr}{tr}
\newcommand{\ntr}{\widetilde{\smash\tr\vphantom{\raisebox{.5pt}{r}}}}
\newcommand{\N}{\mathbb{N}}
\newcommand{\PP}{\mathbb{P}}
\newcommand{\Po}{\mathcal{P}}
\newcommand{\Ca}{\mathfrak{C}}
\newcommand{\id}{\text{id}}
\newcommand{\VV}{\mathbb{V}}
\newcommand{\oVV}{\boldsymbol{\mathaccent'27\VV}}  
\DeclareMathOperator{\Ord}{Ord}
\DeclareMathOperator{\Hom}{Hom}
\newcommand{\sut}{{\mathrm{SU}(2)}}
\newcommand{\Uo}{{\mathrm{U}(1)}}
\newcommand{\Hi}{\mathcal{H}}
\newcommand{\oAA}{\boldsymbol{\mathaccent'27\AA}}  
\newcommand{\Ai}{\A_\infty}
\newcommand{\longinto}{\lhook\joinrel\relbar\joinrel\rightarrow}
\newcommand{\longonfrom}{\twoheadleftarrow\joinrel\relbar}
\newcommand{\XX}[2]{\left<#1\Psi^{N\Lambda}|#2\Psi^{N\Lambda}\right>}
\newcommand{\SO}{\mathrm{SO}}
\DeclareMathOperator{\Dynkin}{Dynkin}

\begin{document}

\notinjournal{
\begin{flushright}
\vspace{-0.5in}
\begin{tabular}{l}
\textsf{\small CGPG-97/8-1}\\
\textsf{\small q-alg/9708030}\\
\end{tabular}
\vspace{0.5in}
\end{flushright}}

\title{Quantization of Equivariant Vector Bundles}
\author{Eli Hawkins}
\notinjournal{}
\maketitle
\begin{center}
\vspace{-4ex}
\small\emph{\small Center for Gravitational Physics and Geometry}\\
\emph{\small The Pennsylvania State University,
University Park, PA 16802}\\
{\small E-mail: mrmuon@phys.psu.edu}\\
\end{center}

\begin{abstract}
The quantization of vector bundles is defined. Examples are 
constructed for the well controlled case of equivariant vector 
bundles over compact coadjoint orbits. (A coadjoint orbit is a 
symplectic manifold with a transitive, semisimple symmetry group.) In 
preparation for the main result, the quantization of coadjoint orbits 
is discussed in detail. 

This subject should not be confused with the quantization of the 
total space of a vector bundle such as the cotangent bundle.
\end{abstract}

\tableofcontents

\section{Introduction}
Quantization is a vaguely defined process by which a noncommutative 
algebra is generated from some ordinary, commutative space. 
Traditionally this space has been the phase space of some system in 
classical mechanics; the algebra is then meant to consist of 
observables for a corresponding quantum system. A more recent use of 
quantization is with a space that is thought of geometrically; the 
quantization is then thought of as giving noncommutative geometries 
which approximate the original space being considered. 

The existing theory of quantization is limited for this purpose in 
that it only gives an algebra. This corresponds to just having the 
topology of the quantized space (see \cite{con1}). If the original 
space has more interesting structures than just its topology, then it 
would be desirable to in some sense ``quantize'' these as well.

Arguably, vector bundles are the most important structures beyond 
point set topology. Most structures used in geometry are, or involve, 
vector bundles. The vector fields, differential forms, and spinor 
fields are sections of vector bundles. $K$-theory is constructed from 
vector bundles. A Riemannian metric is a section of a bundle. 
Differential operators, such as the Dirac operator, act on sections 
of vector bundles. Indeed, in physics most fields are sections of 
vector bundles.

This paper is a first step towards a theory of the 
quantization of vector bundles. In pursuit of this goal, I present a 
plausible definition for the quantization of a vector bundle, and 
illustrate it with a large class of  examples. I give a more 
general construction of quantization of vector bundles in \cite{haw4}.
\gap

I only consider compact manifolds in this paper for several reasons. 
One is that this is inevitably the simplest case to deal with, since 
almost anything that will work generally will work in the compact 
case. Another is physically motivated. The most natural quantizations 
of compact manifolds give finite-dimensional algebras; as a result, 
the degrees of freedom of anything on the space should become finite 
after quantization. This can, therefore, be used as a regularization 
technique for quantum field theories (see \cite{g-k-p1}). Outside of 
some definitions, I will assume the space is a compact manifold $\M$ 
and the quantizations are finite-dimensional.

In order to get simple examples I will assume that the geometry is 
also highly symmetrical. Suppose that some compact, semisimple Lie 
group acts transitively\footnote{I.~e., the group can take any point 
to any other point.} on $\M$, and that everything is equivariant 
under the action of this group. A manifold that can be quantized (to 
give finite-dimensional algebras) in a reasonable sense must have a 
symplectic structure (App.~\ref{symplectic}). A symplectic manifold 
with transitive symmetry by a compact, semisimple Lie group must be 
equivalent to a coadjoint orbit of that group (App.~\ref{coadjoint}). 
The coadjoint orbits are therefore the only spaces that can be 
quantized nicely with this much symmetry. Luckily, coadjoint orbits 
of compact Lie groups have a very simple systematic quantization 
(Sec.~\ref{Q_co_orb}).
\gap

I begin in Section \ref{generalQ} with a general definition of 
quantization structure similar to that given by Berezin in 
\cite{ber}. This definition involves a minimum of structure. However, 
greater structure can be useful for some purposes.

The perspective of noncommutative geometry \cite{con1} holds that a 
noncommutative algebra should correspond to the ``true'' geometry, 
and that the ``classical limit'' is merely a convenient approximation to 
this \cite{c-c1,con2}. This suggests that the classical algebra of 
functions should be secondary, constructed as the limit of a sequence 
of noncommutative algebras. Based on this philosophy (and other 
motivations described in Sec.~\ref{conclusions}), in Section 
\ref{limitQ} I outline an approach to quantization based on a 
directed or inverse system of algebras whose limit is the classical 
algebra of continuous functions; I call these structures direct and inverse 
limit quantizations. The technical details of these limits are 
discussed in Appendices \ref{alg_dir_lim} and \ref{alg_inv_lim}.

In Section \ref{vector}, I give a definition for the quantization of 
a vector bundle. Like the quantization of an algebra, the 
quantization of a vector bundle can be viewed in terms of a directed 
or inverse system. This is described in Section \ref{limit_vector} 
and detailed in Appendices \ref{mod_dir_lim} and \ref{mod_inv_lim}.

The most relevant properties of homogeneous spaces and their vector 
bundles are described in Section \ref{homogeneous}. In Appendix 
\ref{restrictions}, I describe the reasons that the spaces considered 
here are all coadjoint orbits, and then discuss some properties of 
these spaces. Appendix \ref{coadjoint_classification} describes the 
classification of the coadjoint orbits for a given group, and gives a 
diagrammatic technique for expressing a coadjoint orbit as a coset 
space.

The standard quantization of coadjoint orbits is reviewed, and 
described in perhaps new ways, in Section \ref{Q_co_orb}. The 
quantization is constructed using generators and relations in Section 
\ref{gen_pic}. The directed and inverse limit quantizations are 
constructed in \ref{coherent}. Appendix \ref{projective} gives some 
additional details which are relevant to the discussion of 
convergence of the direct and inverse limit quantizations in 
\ref{convergence}.

Section \ref{Q_vec_buns} contains the main results of this paper. I 
first construct quantized vector bundles, and then determine what 
bundles these are quantizations of. I show that all equivariant 
vector bundles over coadjoint orbits may be quantized.

I then discuss some matters arising from this construction. In 
\ref{unique} I explain the extent to which the construction is 
unique. In \ref{geoQ} I note an interesting relationship to geometric 
quantization. In \ref{bimodules} I note a property that these 
quantizations fail to have.

In order to illustrate the constructions in this paper, I describe 
some of the details in the simplest possible case, that of $S^2$, in 
Section \ref{twosphere}.

Appendices \ref{sections} and \ref{representation_theory} serve to 
fill in some background and fix notation. Appendix \ref{sections} is 
background mainly for Appendix \ref{lims}\@. Some of the relevant 
facts about Lie groups are reviewed in Appendix 
\ref{representation_theory} in a perspective appropriate to this 
paper.

This topic unfortunately requires using a great many symbols. A table 
of notations is provided at the end of the paper.

\section{Quantization}
\label{generalQ}
Generally, quantization refers to some sort of correspondence between 
an algebra of functions on some space, and some noncommutative 
algebra. This might involve a map that identifies functions to 
operators in the noncommutative algebra, or perhaps \emph{vice 
versa}. The idea of a ``classical limit'' is that the algebra of 
quantum operators becomes the algebra of classical functions in some 
limiting sense. To make this meaningful requires having not one, but 
a whole sequence (discrete or continuous) of quantum algebras. 

This idea can be made more concrete. Let all algebras involved be 
\cs-algebras. Call the space $\M$; the algebra of functions is the 
algebra $\C_0(\M)$ of continuous functions (vanishing at infinity in 
the noncompact case). The set of quantum algebras may be parameterized 
either continuously (say, over $\I=\Rplus$) or discretely (say, over 
$\I=\N$). Compactify the parameter space $\I$ by adjoining some 
``$\infty$'' where the classical limit belongs. The algebras form a 
bundle $\A_{\hat\I}$ over this completed parameter space 
$\hat\I=\I\cup\{\infty\}$, each quantum algebra is the fiber over its 
parameter and $\C_0(\M)$ is the fiber over $\infty$. This 
$\A_{\hat\I}$ should in fact be a continuous field of \cs-algebras; 
see \cite{dix}.

I am taking the perspective in this paper that quantization gives 
noncommutative approximations to the topology $\M$. From this 
perspective, the most essential information about the 
quantum-classical correspondence is encoded in the topology of the 
bundle $\A_{\hat\I}$. A sequence of operators in each of the quantum 
algebras can be reasonably identified with a certain function only if 
these together form a continuous section of $\A_{\hat\I}$. The space of 
continuous sections over $\hat\I$ is naturally a \cs-algebra, 
$\AA:=\G(\A_{\hat\I})$ (see App.~\ref{sections}). There is a natural 
surjection $\Po:\AA\onto\C_0(\M)$ which is simply evaluation at the 
point $\infty\in\hat\I$. This algebra and surjection are the most 
succinct and bare-bones quantization structure. This will be referred 
to as a \emph{general quantization}. This is almost the same as the 
structure of quantization given by Berezin in \cite{ber}\footnote{The 
major difference is that Berezin used smooth rather than continuous 
functions.}\negthinspace. It is also a generalization of the structure of a 
strict deformation quantization \cite{rie}; in that case the index 
set $\hat\I$ is required to be an interval.

Other quantization structures contain more (possibly irrelevant) 
information. Suppose that we are given a quantization of a space $\M$ 
in the form of a sequence of algebras $\{\A_N\}_{N=1}^{\;\infty}$ and 
maps $P_N:\C_0(\M)\onto\A_N$. This is a pretty typical quantization 
structure; the operator $P_N(f)$ is considered to be the quantization 
of the function $f$. The topology I give to 
$\A_{\hat\N}=\A_\N\cup\C(\M)$ is the weakest such that for each 
$f\in\C_0(\M)$ the section taking $N\mapsto P_N(f)$ and 
$\infty\mapsto f$ is continuous. Two sets of $P_N$'s that give the 
same topology to the bundle are equivalent for the purposes of my 
perspective.

This structure of general quantization is not tied to any particular 
method of quantization. Indeed, it need not correspond to something 
that would usually be called quantization. The point of it is that a 
large class of concepts of quantization can be used to construct a 
general quantization structure, and it is this structure which is 
relevant to defining the quantization of a vector bundle in Section 
\ref{vector}.

The strategy for constructing general quantizations that is used here 
is that $\AA\equiv\G(\A_{\hat\I})$ is a subalgebra of $\Gb(\A_\I)$ 
(the \cs-algebra of bounded sections over $\I$; see 
App.~\ref{sections}). The difference between these two types of 
sections is the behavior approaching $\infty$; elements of 
$\G(\A_{\hat\I})$ must be continuous at $\infty$. The key is to 
describe the condition of continuity at $\infty$ purely in terms of 
$\I\not\ni\infty$.

\subsection{Direct and inverse limit quantization} 
\label{limitQ}
In this section I make the assumption that $\M$ is compact and the 
quantum algebras are finite-dimensional. Since dimensions change 
discretely, the simplest choice of parameter space is $\I=\N$.

One perspective on quantization is that the classical algebra is 
\emph{literally} the limit of the sequence of quantum algebras. A 
limit of algebraic objects is generally constructed from either a 
``directed system'' or ``inverse system''\negthinspace, so those are 
what I use here. The former is a bundle of algebras $\A_\N$ and a sequence 
of maps $i_N:\A_N\into\A_{N+1}$ linking them together. In the latter the 
maps are in the opposite direction, $p_N:\A_N\onto\A_{N-1}$. If 
constructed properly, these types of systems have limits 
$\dlim\{\A_*,i_*\}$ and $\ilim\{\A_*,p_*\}$ which are \cs-algebras; 
these are detailed in Appendices \ref{alg_dir_lim} and 
\ref{alg_inv_lim}. 

Intuitively, the directed system can be thought of as
\beq
\A_1\stackrel{i_1}{\longinto}\A_2\stackrel{i_2}{\longinto} \ldots 
\longinto\dlim\{\A_*,i_*\}
\;.\eeq
For every $N$ there is a composed injection 
$I_N:\A_N\into\dlim\{\A_*,i_*\}$. These satisfy a consistency 
condition with the $i_N$'s that $I_N=I_{N+1}\circ i_N$. Similarly, 
the inverse system can be thought of as 
\beq
\A_1\stackrel{p_2}{\longonfrom}\A_2\stackrel{p_3}{\longonfrom} \ldots 
\longonfrom \ilim\{\A_*,p_*\}
\;.\eeq
There are composed surjections $P_N:\ilim\{\A_*,p_*\}\onto\A_N$. 
These also satisfy a consistency condition that $P_N=p_{N+1}\circ 
P_{N+1}$. These $I_N$'s and $P_N$'s are part of the general 
constructions of directed and inverse limits.
The general quantization algebra $\AA$ is also a natural byproduct of 
these constructions.

The maps $i_N$ and $p_N$ used in these must not be assumed to be 
(multiplicative) homomorphisms in general. That assumption would 
actually restrict $\M$ to be a totally disconnected space, which is 
almost certainly not what we want. Instead we must allow these maps 
to be some more general type of morphisms, such as unital completely 
positive maps\footnote{The property of complete positivity will not 
be used here; although it will be mentioned several times. For 
definition and discussion see \cite{lan}.}\negthinspace. This is 
discussed a little more in Appendix \ref{alg_dir_lim}. 

\section{Quantized Vector Bundles}
\label{vector}
Suppose that we are given a finitely generated vector bundle 
$V\onto\M$ (see \cite{sch}). If the algebra of functions $\C_0(\M)$ 
is quantized, then what should be meant by the quantization of $V$? 
In noncommutative geometry, all geometrical structures are dealt with 
algebraically. In order to find the appropriate definition for 
quantization of $V\!$, we must first treat $V$ algebraically. The 
algebraic approach comes from the fact that the continuous sections 
$\Gn(V)$ form a finitely generated, projective module of the algebra 
$\C_0(\M)$. Indeed, this gives a one-to-one correspondence between 
finitely generated, locally trivial, vector bundles and finitely 
generated, projective modules (see \cite{con1}).
The ``quantization'' of $V$ should give modules for each of the 
quantum algebras $\A_N$; in other words, a bundle of modules over 
$\I$. 

I define a quantization of the bundle $V$ to be a bundle of modules 
$V_{\hat\I}$ over $\hat\I$ such that the topology is consistent with 
that of $\A_{\hat\I}$, and the fiber at $\infty$ is the module 
$\Gn(V)$. 

The space of sections $\VV:=\G(V_{\hat\I})$ is a module of $\AA$. 
This gives another way of describing the quantization of $V\!$. A 
quantization of $V$ may be equivalently defined as a finitely 
generated, projective module $\VV$ of $\AA$ satisfying the sole 
condition that the push-forward by $\Po$ to a module of $\C_0(\M)$ is 
$\Gn(V)$. The condition that $\A_{\hat\I}$ and $V_{\hat\I}$ have 
consistent topologies is implicitly encoded in this definition.

Just as a continuous function is not uniquely determined by its value 
at a single point, there is not a single, unique quantization of a 
given $V\!$. Indeed, when $\I$ is discrete, any finite subset of 
$V_N$'s can be changed arbitrarily. However, there may be a uniquely 
natural choice for almost all $V_N$'s given by a single formula. This 
is so in the case discussed in this paper. This issue is discussed 
further in Section \ref{unique}.
\gap

The guiding principle for quantizing vector bundles will be that we 
already have one example. The sections of the trivial line bundle 
$V=\M\times\co$ are simply the continuous functions $\C_0(\M)$. This 
means that $\VV=\AA$ should always be a good quantization of this 
bundle.

\subsection{Direct and inverse limits}
\label{limit_vector} 
Return to the assumptions of Section \ref{limitQ} (compactness, 
etc.). As with quantizing $\C(\M)$, it is possible to use additional 
structure in the quantization of a vector bundle. A quantized vector 
bundle can be constructed from a directed system $\{V_*,\iota_*\}$ or 
an inverse system $\{V_*,\pi_*\}$ of modules. In these systems, each 
$V_N$ is an $\A_N$-module; the maps are linear maps $\iota_N:V_N\into 
V_{N+1}$ and $\pi_N:V_N\onto V_{N-1}$. The details of this are 
described in Appendices \ref{mod_dir_lim} and \ref{mod_inv_lim}. 
There are again composed injections $I^V_N$ and surjections $P^V_N$, 
satisfying the same sort of compatibility conditions as for $I_N$ and 
$P_N$ in Section \ref{limitQ}.

\section{Classical Homogeneous Spaces}
\label{homogeneous}
Again, and throughout the rest of this paper, I assume that $\M$ is a 
compact manifold, the parameter space is $\I=\N$, and the algebras 
$\A_N$ are finite-dimensional. 
In order to get some control of the system, and construct some 
quantizations explicitly, let us assume that some group $G$ acts 
transitively on $\M$ (i.~e., $\M$ is homogeneous) and that everything 
we do will be $G$-equivariant. It is a standard construction (see 
\cite{hel}) that $\M$ can be written as a coset space $\M=G/H$ 
where the isotropy group is $H:=\{h\in G\mid h(o)=o\}$ for some 
arbitrary basepoint $o\in\M$. Since $\M$ is a manifold, $G$ is best 
chosen to  be a Lie group. If we assume $G$ to be compact and 
semisimple\footnote{Assuming $G$ semisimple is equivalent to assuming 
$\M$ is not a torus or the product of a torus with something 
else.}\negthinspace, then the set of $\M$'s we are interested in is (up to 
equivalence) the set of ``coadjoint orbits'' (see 
App.~\ref{coadjoint}).

\subsection{The set of coadjoint orbits}
\label{intro_coadjoint}
The coadjoint space is $\g^*$, the linear dual of the Lie algebra 
$\g$ of $G$. There is a natural, linear action of $G$ on $\g^*$. A 
coadjoint orbit is simply the orbit of some point in $\g^*$ under 
that $G$ action.

The relevant definitions concerning Lie groups are summarized in 
Appendix \ref{representation_theory}. The classification of coadjoint 
orbits is strikingly similar to the classification of irreducible 
representations. The irreducible representations are classified by 
the dominant weights, which are the vectors on the weight lattice 
that lie in the positive Weyl chamber $\C_+\subset\g^*$. The 
coadjoint orbits are classified by \emph{all} vectors in $\C_+$ (see 
App.~\ref{coadjoint_classification}). Denote by $\Orl$ the coadjoint 
orbit of $\Lambda\in\C_+\subset\g^*$.

Since a coadjoint orbit is a homogeneous space, it can always be 
expressed as a coset space $\Orl\cong G/H$; it is natural to identify 
the basepoint $o=eH\in  G/H$ with $\Lambda\in\Orl$. A diagrammatic 
method of calculating $H$ from $\Lambda$ is described in Appendix 
\ref{coadjoint_classification}. 

The structures of the sets of irreducible representations of $G$ and 
of $H$ are closely related. The weight lattices of $G$ and $H$ are 
naturally identified. However, the sets of weights which are dominant 
(and thus actually correspond to representations) are different. This 
is relevant in Section \ref{identification}.

\subsection{Equivariant bundles}
\label{equi_bun}
\begin{ntn}
In this paper I will generally refer to a representation space (group 
module) simply as a representation. 
\end{ntn}

Suppose that $V$ is an equivariant vector bundle over $\M=G/H$. This 
simply means that $\G(V)$ is a representation of $G$. The fiber $V_o$ 
at the basepoint $o=eH$ is a vector space and is acted on by $H$, so 
$V_o$ is a representation of $H$. 

Suppose that $W$ is a representation of $H$. The set $W\times_HG := 
W\times G/{\sim}$, where $(w,g)\sim(hw,gh^{-1})$, is naturally an 
equivariant vector bundle over $\M$. The bundle surjection 
$W\times_HG\onto G/H$ is $[(w,g)]\mapsto gH$; the action of $g'\in G$ 
is $[(w,g)]\mapsto[(w,g'g)]$. Up to equivalence, all equivariant 
vector bundles may be constructed in this way.

The fiber of $W\times_HG$ at $o$ is simply $W$, so there is a 
one-to-one correspondence between $H$-representations and equivariant 
vector bundles over $\M$.
The semigroup of equivariant vector bundles under direct sum is 
generated by the set of \emph{irreducible} bundles --- those 
corresponding to irreducible representations.

It is not the case that all vector bundles over $\M$ can be made 
equivariant. Nevertheless, I am only considering equivariant bundles 
in this paper. Every bundle over a homogeneous space which is 
mentioned in this paper is a finitely generated, locally trivial, 
equivariant, vector bundle; but I will frequently omit some of these 
adjectives.

\section{Quantized Coadjoint Orbits}
\label{Q_co_orb}
\begin{ntn}
The irreducible representations of $G$ are in one-to-one 
correspondence with dominant weights 
(App.~\ref{representation_theory}). Denote the space of the 
representation corresponding to the weight $\lambda$ by $(\lambda)$. 
This is the $G$-representation with ``highest weight'' $\lambda$ 
(App.~\ref{representation_theory}).

Denote $\A_N:=\End (N\Lambda)$, the algebra of matrices on the vector 
space $(N\Lambda)$; the notation $\A_N$ will be justified in the 
following.
\end{ntn}

\subsection{Generators and relations picture}
\label{gen_pic}
The action of $\g$ on $(N\Lambda)$ can be expressed as a map 
$\g\to\End(N\Lambda)=\A_N$. The associative algebra $\A_N$ is 
generated by the image of the Lie algebra $\g$. Let 
$\{J_i\}\subset\g$ be a basis of self-adjoint generators of $\g$ 
acting on $(N\Lambda)$; $\A_N$ can be written in terms of this set of 
generators and the following relations. 

First, the commutation relations state that
\beq
[J_i,J_j]_-=iC^k_{\;ij}J_k
\label{J_commutation}\;,\eeq
where $C^k_{\;ij}$ are the structure coefficients. Second, the 
Casimir relations state that
\beq
\C_n(J)=c_n(N\Lambda) \quad \forall n
\label{J_Casimir},\eeq
where the Casimirs $\C_n$ are $G$-invariant, symmetrically ordered, 
homogeneous polynomials in the $J$'s, and the $c_n$'s are the 
corresponding eigenvalues. Finally, the Serre relations state that 
certain linear combinations of $J_i$'s are nilpotent, the order of 
nilpotency rising linearly with $N$; an example of this is given in 
Section \ref{twosphere}.

The Casimir eigenvalues $c_n(N\Lambda)$ are polynomials in $N\Lambda$ 
of the same order as $\C_n$. In fact the leading order (in $N$) term 
is $\C_n(\Lambda)N^{\Ord(\C_n)}$. The reason that it is meaningful to 
evaluate $\C_n$ on a point of $\g^*$ (such as $\Lambda$) as well as 
on the $J_i$'s is that the $J_i$'s together form a sort of Lie algebra 
valued vector in $\g^*$.

The Serre relations are actually equivalent to the condition that the 
$J_i$'s generate a \cs-algebra. Suppose that the $J_i$'s do lie 
inside a \cs-algebra and satisfy the commutation and Casimir 
relations. Then this \cs-algebra can be faithfully represented on a 
Hilbert space $\Hi$. The commutation relations imply that the $J_i$'s 
generate a unitary representation of $G$ on $\Hi$. The Casimir 
relations imply that $\Hi$ can only be $(N\Lambda)$ or some Hilbert 
space direct sum of copies of $(N\Lambda)$. This means that the 
\cs-subalgebra generated by the $J_i$'s is $\End(N\Lambda)$; which 
implies that the Serre relations are satisfied.
\gap

Now, regard the $\A_N$'s as forming a bundle $\A_\N$ over the 
discrete parameter space $\N$. We can think of $N$ and the generators 
$J_i$ as sections in $\G(\A_\N)$, but neither is bounded, so they are 
not in $\Gb(\A_\N)$ (the \cs-algebra of bounded sections; see 
App.~\ref{sections}). However, the combinations $X_i=N^{-1}J_i$ are 
bounded; as can be seen by considering the quadratic 
Casimir\footnote{The eigenvalue of the quadratic $\C_1(X)$ is 
$\C_1(\Lambda)$ plus a term proportional to $N^{-1}$, therefore it is 
bounded as $N\to\infty$, therefore it is a polynomial of bounded 
operators.} $\C_1$. This means that $X_i\in\Gb(\A_\N)$. 

Define $\AA$ to be the \cs-subalgebra of $\Gb(\AA_\N)$ generated by 
the $X_i$'s. Define $\AA_0:=\Gn(\A_\N)$ to be the algebra of 
sections vanishing at $\infty$  (see App.~\ref{sections}); 
since in fact\footnote{It is essentially sufficient to show that 
$\AA$ contains one function on $\N$ that nontrivially converges to 
$0$.} $\AA_0$ is contained in $\AA$, it is an ideal 
there\footnote{Because $\AA_0$ is an ideal in 
$\Gb(\AA_\N)$.}\negthinspace. 
Define 
$$
\Po:\AA\onto\Ai:=\AA/\AA_0
$$
to be the corresponding quotient homomorphism; this essentially just 
evaluates the $N\to\infty$ limit.

By construction, the images $x_i:=\Po(X_i)$ generate the quotient 
algebra $\Ai$. The relations these satisfy all derive from the 
relations satisfied by the $X_i$'s. These generators commute, since 
\beq
[x_i,x_j]_-=\Po([X_i,X_j]_-)=\Po(iN^{-1}C^k_{\;ij}X_k)=0
\;,\eeq
so $\Ai$ is a commutative \cs-algebra (and therefore is the algebra 
of continuous functions on some space). The $x_i$'s transform under 
$G$ in the same way as Cartesian coordinates on $\g^*\!$, so $\Ai$ is 
the algebra of continuous functions on some subspace of $\g^*\!$. The 
non-Serre relations alone define a \cs-algebra; therefore the Serre 
relations do not give any additional relations for $\Ai$. The only 
other relations the $x_i$'s satisfy are polynomial relations
\beq
\C_n(x) \stackrel{!}= \lim_{N\to\infty} N^{-\Ord(\C_n)} c_n(N\Lambda) 
=\C_n(\Lambda)
\eeq
which make $\Ai$ the algebra of continuous functions on the algebraic 
subspace $\M\subset\g^*$ determined by these polynomials. 

The Casimir polynomials are a complete system of $G$-invariant 
polynomials; therefore $\M$ must be a single coadjoint orbit. 
Obviously, $x=\Lambda$ satisfies $\C_n(x)=\C_n(\Lambda)$, so 
$\Lambda\in\M$; therefore $\M$ is the orbit $\Orl$.  This shows that 
in the sense of Section \ref{generalQ}, the system 
$\Po:\AA\onto\C(\Orl)$ is a general quantization of $\Orl$.
\gap

In this construction the $\Lambda$ was required to be integral (a 
weight) rather than any arbitrary  $\Lambda\in\C_+$. However, this is 
not a serious restriction. Rescaling $\Lambda$ simply rescales 
$\Orl$, therefore a more appropriate parameter space for distinct 
coadjoint orbits is the projectivisation $\PP\C_+$. The image of the 
weights is dense in $\PP\C_+$ (it is the set of ``rational'' points), 
so the quantizable coadjoint orbits are dense in the space of 
distinct coadjoint orbits.

\subsection{Limit quantization picture}
\label{coherent}
\begin{ntn}
The linear dual of an irreducible representation is also an 
irreducible representation; we can therefore define $\lambda^*$ by 
the property $(\lambda^*)=(\lambda)^*$. This is a linear 
transformation on the weights (see App.~\ref{representation_theory}). 
With this notation $\A_N\equiv\End(N\Lambda) = 
(N\Lambda)\otimes(N\Lambda^*)$.

Given a choice of Cartan subalgebra and positive Weyl chamber, there 
is a preferred, 1-dimensional ``highest weight subspace'' in 
$(N\Lambda)$; choose a normalized basis vector $\Psi^{N\Lambda}$ 
there and call it the highest weight vector (see 
App.~\ref{representation_theory}). 
\end{ntn}

Not only do the coadjoint orbits have equivariant general 
quantizations, but they also admit equivariant direct and inverse 
limit quantizations. There are standard constructions of maps 
$\A_N\into\C(\Orl)$ and $\C(\Orl)\onto\A_N$ which are suitable to be 
used as $I_N$ and $P_N$. I present these first.
\gap

We need an equivariant, linear injection $I_N:\A_N\into\C(\Orl)$. If 
we have such an $I_N$, then for every point $x\in\Orl$, evaluation at 
$x$ determines a linear function 
$$
I_N(\,\cdot\,)(x):\A_N\to\co
\;;$$ 
in other words, $x$ gives an element of the dual $\A_N^*$. Such an 
$I_N$ is in fact equivalent to an injection 
$I^*_N:\Orl\into\A_N^*=(N\Lambda^*)\otimes(N\Lambda)$. Since $I^*_N$ 
must be equivariant, it is completely specified by the image of the 
basepoint $o=eH$. This image must be $H$-invariant.

The highest weight vector $\Psi^{N\Lambda}\in(N\Lambda)$ is 
$H$-invariant, modulo phase. Its conjugate vector 
$\Psi^{-N\Lambda}\in(N\Lambda^*)$ transforms by the opposite phase, 
so the product $\Psi^{-N\Lambda}\otimes\Psi^{N\Lambda}\in\g^*$ is 
$H$-invariant. In fact, $H$ is the largest subgroup that this is 
invariant under. Define the image of the basepoint to be 
$I^*_N(o):=\Psi^{-N\Lambda}\otimes\Psi^{N\Lambda}\in\g^*$. With this 
choice, $I_N$ is given by 
\beq
I_N(a)(gH)=\left< g\Psi^{N\Lambda}\right| a\left| 
g\Psi^{N\Lambda}\right>
\label{I_map}\eeq
 for any $gH\in\Orl$.

There is some apparent arbitrariness in this construction. There were 
choices made of Cartan subalgebra, positive Weyl chamber, and phase 
of the highest weight vector. However, the resulting $I_N$ is only 
arbitrary by the freedom to rotate $\Orl$ about $o$ (by $H$), and 
this freedom was inevitable. 

We now need to construct injections $i_N:\A_N\into\A_{N+1}$. The 
question is how to get from something acting on $(N\Lambda)$ to 
something acting on $([N+1]\Lambda)$. The key is that precisely one 
copy of $([N+1]\Lambda)$ always occurs as a subrepresentation of 
$(\Lambda)\otimes(N\Lambda)$ (see App.~\ref{representation_theory}). 
There is a unique, natural projection 
$$
\Pi_+\in\Hom_G\left[(\Lambda)\otimes(N\Lambda),\left([N+1]\Lambda\right)\right]
$$
which maps a vector in $(\Lambda)\otimes(N\Lambda)$ to its component 
in the irreducible subrepresentation 
$([N+1]\Lambda)\subset(\Lambda)\otimes(N\Lambda)$. Using $\Pi_+$, an 
element $A\in\End[(\Lambda)\otimes(N\Lambda)]$ can be mapped to 
$\Pi_+ A \,\Pi_+^* \in \A_{N+1}$. Now, that algebra is 
$$
\End[(\Lambda)\otimes(N\Lambda)]= \End(\Lambda)\otimes\End(N\Lambda) 
= \A_1\otimes\A_N
\;.$$ 
There is a very simple map $\A_N\into\A_1\otimes\A_N$ taking 
$a\mapsto 1\otimes a$. Composing these gives, as desired, a map 
$i_N:\A_N\into\A_{N+1}$ by the formula 
\beq
i_N(a)=\Pi_+(1\otimes a)\Pi_+^*
\label{alg_i}\;.\eeq
This (and any map that can be written in this form) is a completely 
positive map (see \cite{lan}). 

To verify that our $i_N$ really satisfies the consistency condition 
$I_{N+1}\circ i_N = I_N$, it is sufficient to check this at the 
basepoint $o\in\Orl$. So, $\forall a\in\A_N$
\begin{align*}
I_{N+1}\circ i_N(a)(o) 
&= \left< \Psi^{(N+1)\Lambda}\right| i_N(a) 
\left|\Psi^{(N+1)\Lambda}\right> 
\\&= \left< \Psi^{(N+1)\Lambda}\right| \Pi_+(1\otimes a)\Pi_+^* 
\left|\Psi^{(N+1)\Lambda}\right> 
\\&= \left< \Psi^1\otimes\Psi^{N\Lambda}\right| (1\otimes a) 
\left|\Psi^1\otimes \Psi^{N\Lambda}\right>
\\&= \left< \Psi^{N\Lambda}\right| a \left|\Psi^{N\Lambda}\right>
= I_N(a)(o)
\end{align*}
and it is consistent.
\gap

The surjections come about similarly. There is a related function 
$e_N$  taking $\Orl$ to projections in $\A_N$. This maps 
$e_N:o\mapsto \lvert\Psi^{N\Lambda}\rangle \langle 
\Psi^{N\Lambda}\rvert$. Using $e_N$, the injection $I_N$ can be 
written as 
\beq
I_N(a)(x) = \tr [a\,e_N(x)]
\label{pformI}\;;\eeq
and the surjection $P_N$ is defined as
\beq
P_N(f) = \dim(N\Lambda) \int_{\Orl} fe_N\epsilon
\label{pformP}\;,\eeq
where $\epsilon$ is an invariant volume form normalized to give 
$\Orl$ volume 1. This map is unital and positive. It is actually the 
adjoint of the map $I_N$ if we put natural inner products on $\A_N$ 
and $\C(\Orl)$. The inner product on $\A_N$ is $\langle a,b\rangle = 
\ntr_{(N\Lambda)}(a^*b)$; where $\ntr_{(N\Lambda)}$ is the trace over 
$(N\Lambda)$, normalized to give $\ntr_{(N\Lambda)} 1 = 1$. The inner 
product on $\C(\M)$ is $\langle f_1,f_2\rangle = \int_{\Orl} 
f_1^*f_2\epsilon$.

We will automatically satisfy the consistency with the  $P_N$'s if we 
choose $p_N$ to be the adjoint of $i_{N-1}$. The immediately obtained 
formula is 
\begin{subequations}
\beq
p_N(a)=[\ntr_{(\Lambda)}\otimes \id_{\A_{N-1}}](a\oplus 0)
\label{trace_form}\;;\eeq
where this is a partial trace of the action of $a$ on 
$(N\Lambda)\subset(\Lambda)\otimes([N-1]\Lambda)$. This can actually 
be written in essentially the same form as the $i_N$'s.
Precisely one copy of $([N-1]\Lambda)$ always occurs as a 
subrepresentation of $(\Lambda^*)\otimes(N\Lambda)$, so there is a 
corresponding projection 
$\Pi_-\in\Hom_G[(\Lambda^*)\otimes(N\Lambda),([N-1]\Lambda)]$. With 
this, define $p_N:\A_N\onto\A_{N-1}$ by 
\beq
p_N(a) = \Pi_-(1\otimes a)\Pi_-^*
\;.\eeq
\end{subequations}
To see that this is equivalent to \eqref{trace_form}, it is 
sufficient to check that these agree for 
$a=e_N(o)=\left|\Psi^{N\Lambda}\right>\left<\Psi^{N\Lambda}\right|$. 
These $p_N$'s are also completely positive.

\subsection{Convergence}
\label{convergence}
I will now show that these direct and inverse limit quantizations are 
both convergent by considering the ``product'' 
$I_N[P_N(f_1)P_N(f_2)]$ for any two functions $f_1,f_2\in\C(\Orl)$. 
This is not an associative product (compare Eq.\ \eqref{star}), since 
$P_N\circ I_N \neq\id$, but as $N\to\infty$ it nevertheless converges 
to the product of functions.

This ``product'' can be written in terms of an integration kernel as 
\beq
I_N[P_N(f_1)P_N(f_2)](x) = \int\!\!\!\int_{\Orl} K_N(x,y,z) f_1(y) 
f_2(z) \epsilon_y \epsilon_z
\;.\eeq
The volume form $\epsilon$ is again the $G$-invariant volume form 
giving $\Orl$ total volume $1$. From the construction of the maps 
$I_N$ and $P_N$ in \eqref{pformI} and \eqref{pformP} it is immediate 
that 
\beq
K_N(x,y,z) = [\dim(N\Lambda)]^2 \tr [e_N(x)e_N(y)e_N(z)]
\;.\eeq
If we use the identification $\Orl=G/H$, this can be factorized as
\beq
K_N(gH,g'H,g''H) = [\dim(N\Lambda)]^2 
\XX{g}{g'}\XX{g'}{g''}\XX{g''}{g}
\label{kernel}\;.\eeq
The factor of $[\dim (N\Lambda)]^2$ serves to normalize $K_N$ so that 
\beq
I_N[P_N(1)P_N(1)]=1
\label{normalization}\;,\eeq 
as it should be since $P_N(1)=1$ and $I_N(1)=1$. 

The inner products in \eqref{kernel} have several nice properties. By 
construction, these are certainly smooth functions. The absolute 
value $\left|\XX{g}{g'}\right|$ only depends on the points 
$gH,g'H\in\Orl$, and is equal to $1$ for $gH=g'H$; but for any 
$gH\neq g'H$, 
$$
\left|\XX{g}{g'}\right|<1
\;.$$
The fact that (see App.~\ref{representation_theory}) 
$\Psi^{N\Lambda}=\Psi^\Lambda\otimes\dots\otimes\Psi^\Lambda$, gives 
the convenient identity 
\beq
\XX{g}{g'} = \left[\left< g\Psi^\Lambda|g'\Psi^\Lambda\right>\right]^N
\label{convenient}\;.\eeq

These properties imply that for any $gH\neq g'H$, 
$$
\XX{g}{g'}\toas 0
$$
exponentially. The factor $[\dim (N\Lambda)]^2$ only increases 
polynomially; therefore, outside any neighborhood of $x=y=z$, 
$K_N(x,y,z)$ vanishes uniformly as $N\to\infty$. This means that in 
order to investigate the $N\to\infty$ limit, it is sufficient to 
consider $x$, $y$, and $z$ close together.

Since $\Orl$ is homogeneous, we can let $x=o$ without loss of 
generality. In order to construct an approximation for $K_N$ near 
$o$, we need a coordinate patch about $o$. Coadjoint orbits are 
always K\"ahler manifolds, so complex coordinates are convenient. The 
(real) tangent fiber $T_o\Orl$ is naturally a complex Hermitean space 
and in fact can be identified to a subspace of $(\Lambda)$ which is 
orthogonal to $\Psi^\Lambda$. A suitable complex coordinate patch can 
be constructed by using this identification along with the 
exponential map; thus a neighborhood of $o$ is coordinatised by 
vectors in a subspace of $(\Lambda)$. Let $\upsilon$ and $\zeta$ be 
the complex coordinates of $y$ and $z$ respectively. Using these 
coordinates, to second order 
$$
[\dim (\Lambda)]^{-2} K_1(o,y,z) \approx 1 - \norm{\upsilon}^2 - 
\norm{\zeta}^2 + \langle\upsilon|\zeta\rangle
\eqno\eqref{K_approx}$$
(see App.~\ref{projective}). A formula for $K_N$ (with $N\gg1$) can 
be constructed by raising this to the $N$'th power and recalling the 
normalization \eqref{normalization}. This gives
\beq
K_N(o,y,z)\epsilon_y\epsilon_z \approx 
\left(\tfrac{N}{\pi}\right)^{2n} 
e^{-N\left[\lVert\upsilon\rVert^2 + \lVert\zeta\rVert^2 - 
\langle\upsilon|\zeta\rangle\right]} \, d^{2n}\upsilon \, d^{2n}\zeta
\label{Kapprox}\eeq
where $2n = \dim \Orl$.
The $L^1$ norm of the error in this expression is of order 
$N^{-\frac32}$ and thus goes to $0$ as $N\to\infty$. It is a standard 
result that as $N\to\infty$ a complex  Gaussian such as 
\eqref{Kapprox} converges as a $\C^{-\infty}$ distribution to the 
delta distribution $\delta^{2n}(\upsilon)\delta^{2n}(\zeta) \, 
d^{2n}\upsilon \, d^{2n}\zeta$. This means for \emph{smooth} 
functions $f_i\in\C^\infty(\Orl)$ that  $I_N[P_N(f_1)P_N(f_2)](o)\to 
f_1(o)f_2(o)$ as $N\to\infty$, and (using the homogeneity of $\Orl$)
\beq
I_N[P_N(f_1)P_N(f_2)]\toas f_1 f_2
\label{conv}\;.\eeq

If, instead of smooth functions, we have \emph{continuous} functions 
$f_i\in\C(\Orl)$ then we can approximate these with smooth functions 
$\tilde f_i$. Because the maps $I_N$ and $P_N$ are completely 
positive, they are norm-contracting; this implies that the 
norm-difference 
$$
\lVert I_N[P_N(f_1)P_N(f_2)] - I_N[P_N(\tilde f_1)P_N(\tilde f_2)] 
\rVert
$$ 
is bounded uniformly as $N\to\infty$ and goes to $0$ as $\tilde f_i 
\to f_i$. This means that \eqref{conv} is true for all continuous 
functions.

Using the fact that $P_N(1)=1$, this also shows that $I_N$ and $P_N$ 
are asymptotically inverse, in the sense that $I_N\circ P_N(f) \to f$ 
as $N\to\infty$. 
This property means that we can replace $P_N$ by a left inverse of 
$I_N$, and Eq.~\eqref{conv} will continue to hold. This shows 
that the direct limit converges (see App.~\ref{alg_dir_lim}). 
Likewise, we can replace $I_N$ by a right inverse of $P_N$, and 
Eq.~\eqref{conv} will continue to hold. This shows that the 
inverse limit converges (see App.~\ref{alg_inv_lim}).

\subsection{Polynomials}
\label{polynomials}
In Appendix \ref{alg_dir_lim}, the limit $\dlim\{\A_*,i_*\}$ is 
constructed by first constructing the limit $\veclim\{\A_*,i_*\}$ as 
a sequence of vector spaces and then completing to a \cs-algebra. In the 
particular case of coadjoint orbits, $\veclim\{\A_*,i_*\}$ is itself 
interesting. 

The algebra $\C(\Orl)$ is, as a $G$-representation, a closure of the 
direct sum of all its irreducible subrepresentations. On the other 
hand, each $\A_N$ is finite-dimensional and is therefore just a 
direct sum of irreducibles; any element of the limit 
$\veclim\{\A_*,i_*\}$ is in the image of some $\A_N$; therefore, 
$\veclim\{\A_*,i_*\}$ is the ``algebraic'' direct sum of 
irreducibles. Since $\C(\Orl)$ is a closure of this, 
$\veclim\{\A_*,i_*\}$ must be the direct sum of all the irreducible 
subrepresentations of $\C(\Orl)$. 

The polynomial functions $\co[\Orl]$ on $\Orl$ are defined as the 
restrictions to $\Orl$ of polynomials on $\g^*$. The space of 
polynomials of a given degree is a direct sum of irreducible 
representations. Any polynomial has finite degree; therefore 
$\co[\Orl]$ is a direct sum of irreducible representations. Since 
$\co[\Orl]$ is dense in $\C(\Orl)$, it must be the direct sum of the 
irreducible subrepresentations of $\C(\Orl)$. 

This shows that $\veclim\{\A_*,i_*\}=\co[\Orl]$, and so the vector 
space direct limit is in this case an algebra. Whether this is true 
in any more general case remains to be seen.

\section{Quantization of Vector Bundles over $\Orl$}
\label{Q_vec_buns}
\subsection{General quantized bundles}
\label{modules}
$\A_N\equiv\End(N\Lambda)$ is a full (a.~k.~a.\ simple) matrix algebra. The 
classification of the modules of a full matrix algebra is elementary. Any 
module is a tensor product of the fundamental module with 
some vector space. In this case the fundamental module is 
$(N\Lambda)$, and  the vector space should be a $G$-representation. 
Any irreducible, equivariant module of $\A_N$ must be of the form 
\beq
V_N=(N\Lambda)\otimes(\nu)
\;,\eeq
with the algebra only acting on the first factor. Any finitely 
generated, equivariant $\A_N$-module is a direct sum of such 
irreducibles. Because $\A_N$ is finite-dimensional, this $V_N$ is 
automatically projective.

The defining property of a finitely generated, projective module is 
that it is a (complemented) submodule of the algebra $\A_N$ tensored 
with some vector space. This submodule can be picked out by a 
projection (idempotent). In the $G$-equivariant case, ``vector 
space'' becomes ``$G$-representation''\negthinspace, and the projection 
must be $G$-invariant. In the case of this $V_N$, the representation we 
tensor with can be chosen to be irreducible; call it $(\mu)$. This 
means that we can identify $V_N$ with a submodule of 
$\A_N\otimes(\mu)$ in the form
\beq
V_N=[\A_N\otimes(\mu)] \cdot Q_N
\;,\eeq
where $Q_N=Q_N^2$.
The factor $(N\Lambda)$ is treated as a space of column vectors, but 
the factor $(N\Lambda^*)\otimes(\mu)$ is treated as a space of row 
vectors, i.~e., $Q_N$ multiplies them from the right. Acting from the 
left, $Q_N$ would multiply the corresponding (dual) space of column 
vectors $(N\Lambda)\otimes(\mu^*)$; therefore 
$Q_N\in\End[(N\Lambda)\otimes(\mu^*)]=\A_N\otimes\End(\mu^*)$. We can 
choose $\mu$ such that 
$Q_N$ is the unique invariant projection from 
$(N\Lambda)\otimes(\mu^*)$ to the irreducible subrepresentation 
$(\nu^*)$.

The injection $i_N:\A_N\into\A_{N+1}$ can be applied to the tensor 
product of $\A_N$ with a fixed algebra --- in this case 
$\End(\mu^*)$. Let us apply this to $Q_N$ and call the result 
$Q_{N+1}$; by Eq.~\eqref{alg_i}, this is 
\beq
Q_{N+1}:=[i_N\otimes\id](Q_N)= (\Pi_+\otimes1)(1\otimes 
Q_N)(\Pi_+^*\otimes1)
\label{iPiN}\;.\eeq
$Q_{N+1}$ is an endomorphism on $([N+1]\Lambda)\otimes(\mu^*)$ and is 
clearly self-adjoint. Let $\psi\in([N+1]\Lambda)\otimes(\mu^*)$ be a 
normalized vector, and look at the product 
\beq
\langle \psi\rvert Q_{N+1} \lvert\psi\rangle = 
\left<(\Pi_+^*\otimes1)\psi\right| (1\otimes Q_N) 
\left|(\Pi_+^*\otimes1)\psi\right>
\label{Qexp}\;.\eeq
Note that $\Pi_+^*\otimes1$ is just the natural isometric inclusion 
of $([N+1]\Lambda)$ into $(\Lambda)\otimes(N\Lambda)$. The product 
\eqref{Qexp} is equal to $1$ if and only if $(\Pi_+^*\otimes1)\psi$ 
is in the image $(\Lambda)\otimes(\nu^*)$ of $Q_N$; but since 
$\psi\in([N+1]\Lambda)\otimes(\mu^*)$, this is equivalent to $\psi$ 
lying in the intersection $(\Lambda+\nu^*)$. Conversely, \eqref{Qexp} 
is 0 if $(\Pi_+^*\otimes1)\psi$ is orthogonal to 
$(\Lambda)\otimes(\nu^*)$, or equivalently, if $\psi$ is orthogonal 
to $(\Lambda+\nu^*)$. This shows that $Q_{N+1}$ is the projection 
with image $(\Lambda+\nu^*)$.

Note that $\Pi_+^*\Pi_+$ is the self-adjoint idempotent acting on 
$(\Lambda)\otimes(N\Lambda)$, with image $([N+1]\Lambda)$. Using the 
same sort of reasoning as in the last paragraph, the image of 
$(1\otimes Q_N)(\Pi^*_+\otimes1)$ is in 
$([N+1]\Lambda)\otimes(\mu^*)$, so there is the identity
\begin{align}
(1\otimes Q_N)(\Pi^*_+\otimes1) 
&= (\Pi^*_+\Pi_+\otimes1)(1\otimes Q_N)(\Pi^*_+\otimes1) 
\nonumber\\
&= (\Pi^*_+\otimes1)Q_{N+1}
\label{commute_projections}\;.\end{align}
In words, moving $1\otimes Q_N$ right past $\Pi_+^*\otimes1$ 
transforms it into $Q_{N+1}$.

The new projection $Q_{N+1}$ gives an $\A_{N+1}$-module 
$$
V_{N+1}=[\A_{N+1}\otimes(\mu^*)]\cdot Q_{N+1} = 
([N+1]\Lambda)\otimes(\nu+\Lambda^*)
\;.$$ 
Repeating this process gives a whole sequence of modules. Since the 
weight in the second factor is changed by $\Lambda^*$ with each step, 
it is simpler to write in terms of  $\lambda=\nu-N\Lambda^*$. The 
sequence of modules is now 
\beq
V^\lambda_N:=(N\Lambda)\otimes(N\Lambda^*+\lambda)
\label{Vl_form}\;.\eeq
Each of these can be realized as a submodule of $\A_N\otimes(\mu)$ in 
the form
\beq
V^\lambda_N=[\A_N\otimes(\mu)]\cdot Q^\lambda_N
\;.\eeq
The projections are related by the recursion\footnote{Actually, this 
is not \emph{quite} always true; see Sec.~\ref{allow_lambda}.} 
\begin{subequations}\beq
Q^\lambda_{N+1}=[i_N\otimes\id](Q^\lambda_N)
\label{I_Q}\;.\eeq

Because the construction of the $p_N$'s is so similar to that of the 
$i_N$'s, the same reasoning shows that $[p_N\otimes\id](Q_N)$ is a 
projection as well. In fact, the same sequence of projections given 
by \eqref{I_Q} also satisfies 
\beq
Q^\lambda_{N-1}=[p_N\otimes\id](Q^\lambda_N)
\label{P_Q}\;.\eeq\label{I_P_Q}\end{subequations}

Now, we can put all these $Q^\lambda_N$'s together to form 
$Q^\lambda\in \G(\A_\N)\otimes \End(\mu^*)$. The constructions of 
$\AA$ in Appendices \ref{alg_dir_lim} and \ref{alg_inv_lim} say 
essentially that $Q^\lambda\in\AA\otimes\End(\mu^*)$ if and only if 
one of the relations \eqref{I_P_Q} is true in a limiting sense as 
$N\to\infty$. Since Eq.'s~\eqref{I_P_Q} are true for finite $N$, 
we have more than we need to show that 
$Q^\lambda\in\AA\otimes\End(\mu^*)$. 

By construction, this $Q^\lambda$ is obviously a projection. Using 
this, we define 
\beq
\VV^\lambda:=[\AA\otimes(\mu)]\cdot Q^\lambda
\;.\eeq 
This is a well defined, finitely generated, projective module of 
$\AA$, and the restriction to each $\A_N$ is $V^\lambda_N$. This 
shows that $\VV^\lambda$ is a general quantization of \emph{some} 
bundle $V^\lambda$ over $\Orl$. Although $(\mu)$ was used in this 
construction, $\lambda$ completely determines $\VV^\lambda$ as an 
$\AA$-module.

\subsection{Limit quantized bundles}
\label{limit_bun}
 We can use $i_N$ to map 
$$
i_N\otimes\id:\A_N\otimes(\mu)\into\A_{N+1}\otimes(\mu)
\;.$$ 
For some $\psi\in\A_N\otimes(\mu)$, look at what happens to the 
product $\psi Q^\lambda_N$;  using \eqref{commute_projections},
\begin{subequations}\begin{align}
[i_N\otimes\id](\psi Q^\lambda_N) 
&= \Pi_+ (1\otimes[\psi Q^\lambda_N])(\Pi_+^*\otimes 1)
\nonumber\\&
= \Pi_+ (1\otimes\psi)(1\otimes Q^\lambda_N)(\Pi_+^*\otimes 1)
\nonumber\\&
= \Pi_+ (1\otimes\psi)(\Pi_+^*\otimes 1)\cdot 
Q^\lambda_{N+1}\label{handy_form}
\\&= [i_N \otimes\id](\psi)\cdot Q^\lambda_{N+1}
\;.\end{align}\end{subequations}
This implies that $i_N\otimes\id$ maps the image $V^\lambda_N$ of 
$Q^\lambda_N$ to the image $V^\lambda_{N+1}$ of $Q^\lambda_{N+1}$, so 
we can restrict $i_N\otimes\id$ to $V^\lambda_N$ and get a well 
defined injection $\iota_N:V^\lambda_N\into V^\lambda_{N+1}$. These 
injections make a directed system out of the $V^\lambda_N$'s. Because 
of the simple relationship with the directed system $\{\A_*,i_*\}$, 
the system $\{V^\lambda_*,\iota_*\}$ inherits its convergence.

In an essentially identical way, we can construct 
$\pi_N:V^\lambda_N\onto V^\lambda_{N-1}$ as the restriction of 
$p_N\otimes\id$. This gives a convergent inverse system 
$\{V^\lambda_*,\pi_*\}$.

In spite of the way that they were constructed, these $\iota_N$'s and 
$\pi_N$'s are independent of the $(\mu)$ that we use. We can use the 
unique natural projection 
$$
\Pi_{\lambda\,+}\in\Hom_G[(\Lambda)\otimes(N\Lambda+\lambda^*),
([N+1]\Lambda+\lambda^*)]
$$ 
to write (in a slight modification of \eqref{handy_form})
\beq
\iota_N(\psi) = \Pi_+(1\otimes\psi)\Pi_{\lambda\,+}^*
\label{mod_i}\;.\eeq
In this form $\iota_N$ manifestly depends only on $\Lambda$, $N$, and 
$\lambda$. There is again a precisely analogous form for $\pi_N$.

It is easy to see that $\Pi_{0\,+}=\Pi_+$, so $\iota_N$ in 
\eqref{mod_i} is a simple generalization of $i_N$ in \eqref{alg_i}. 
Analogous to the maps $I_N$ and $P_N$ for the algebras, there are 
maps $I^{V^\lambda}_N:V^\lambda_N\into 
V^\lambda_\infty\equiv\G(V^\lambda)$ and 
$P^{V^\lambda}_N:V^\lambda_\infty\onto V^\lambda_N$. These are easily 
constructed as restrictions of $I_N\otimes\id$ and $P_N\otimes\id$.

These limit quantizations both produce the same $\VV^\lambda$ as was 
constructed using $Q^\lambda$ in the previous section. These are, 
therefore, all quantizations of the same bundle $V^\lambda$.

\subsection{Identification with bundles}
\label{identification}
\begin{ntn}
Since the Lie algebras $\g$ and $\h$ share the same Cartan subalgebra, 
their weights are naturally identified (App.~\ref{coadjoint_classification}). 
Denote the $H$-representation with highest weight $\lambda$ by 
$[\lambda]$. Note that $\Psi^\lambda\in[\lambda]\subset(\lambda)$. 
Beware that $[\lambda]^*$ and $[\lambda^*]$ are not generally the 
same.
\end{ntn}

I have established that the irreducible equivariant bundles are given 
by dominant weights of $H$, and irreducible equivariant quantized 
bundles are given by weights of $G$. So, what is the correspondence?

Using the quotient homomorphism $\Po:\AA\onto\C(\Orl)$, define the limit 
projection 
$$
Q^\lambda_\infty:=[\Po\otimes\id](Q^\lambda)\in\C(\Orl)\otimes\End(\mu^*)
\;;$$
this is naturally thought of as a projection-valued function on 
$\Orl$. The bundle $V^\lambda$ can be realized as the subbundle of 
$\Orl\times(\mu)$ determined by $Q^\lambda_\infty$. At each point 
$x\in\Orl$, the fiber of $V^\lambda$ is $V^\lambda_x=(\mu)\cdot 
Q^\lambda_\infty(x)\subset(\mu)$.

The injection $I_N$ is heuristically the limit of applying $i_N$, 
then $i_{N+1}$, then $i_{N+2}$, and so on. The recursion relation 
\eqref{I_Q} thus implies that 
$[I_N\otimes\id](Q^\lambda_N)=Q^\lambda_\infty$.

As explained in Section \ref{equi_bun}, the equivariant bundle 
$V^\lambda$ is completely determined by its fiber at $o\in\Orl$. This 
fiber is given by $Q^\lambda_\infty$ as $V^\lambda_o=(\mu)\cdot 
Q^\lambda_\infty(o)$. It is more convenient to first determine the 
dual $(V^\lambda_o)^*=Q^\lambda_\infty(o)\cdot (\mu^*)$.

The $H$-representation $(V^\lambda_o)^*$ is the image of 
$Q^\lambda_\infty(o)$. This is actually an irreducible 
representation, so it is determined by its highest weight. Let 
$\psi\in(\mu^*)$ be a normalized vector of a given weight. If (and 
only if) $\psi\in(V^\lambda_o)^*$ then 
$\langle\psi|Q^\lambda_\infty(o)|\psi\rangle=1$. So, evaluate this 
expression; it is (using \eqref{I_map})
\begin{align*}
\langle\psi|Q^\lambda_\infty(o)|\psi\rangle 
&= \bigl<\psi\bigr|[(I_N\otimes\id)(Q^\lambda_N)](o)\bigl|\psi\bigr> 
\\&= \left<\Psi^{N\Lambda}\otimes\psi\right| Q^\lambda_N 
\left|\Psi^{N\Lambda}\otimes\psi\right>
\;.\end{align*}
This is 1 if and only if 
$\Psi^{N\Lambda}\otimes\psi\in(N\Lambda+\lambda^*)$. Since 
$N\Lambda+\lambda^*$ is the highest weight of $(N\Lambda+\lambda^*)$, 
the highest weight that $\psi$ can have under this condition is 
$\lambda^*$. This means that $(V^\lambda_o)^*=[\lambda^*]$; therefore 
$V^\lambda_o=[\lambda^*]^*$. Finally, this gives 
\beq
V^\lambda= [\lambda^*]^* \times_H G
\label{the_point}\;.\eeq

\subsection{The allowed weights}
\label{allow_lambda}
The recursion relation \eqref{I_Q} is actually not true for quite all 
values of $\lambda$ and $N$.

If a weight $\nu$ is not dominant, then there really is no 
representation $(\nu)$. It is, however, convenient to define 
$(\nu):=0$ in that case. The condition that $\VV^\lambda\neq0$ is 
that $N\Lambda^*+\lambda$ is dominant for some $N$. 

If $\lambda$ satisfies this condition but is not itself dominant, 
then for low $N$ values $V^\lambda_N=0$, but for sufficiently large 
$N$ values $V^\lambda_N\neq0$. In this case there is some $N$ such 
that $V^\lambda_N=0\neq V^\lambda_{N+1}$. This means that 
$Q^\lambda_N=0\neq Q^\lambda_{N+1}$, so obviously 
$[i_N\otimes\id](Q_N)=0\neq Q^\lambda_{N+1}$ and \eqref{I_Q} fails. 
However, this is the only time that \eqref{I_Q} is not true, so there 
is no real trouble from this. Equation \eqref{P_Q}, on the other 
hand, is always true.

The condition that $V^\lambda$, as given by \eqref{the_point}, is a 
nonzero bundle is that $\lambda^*$ is dominant as an $H$-weight. This 
is actually exactly equivalent to the condition just described for 
$\VV^\lambda\neq0$. 
\emph{This means that any finitely generated, locally trivial, 
equivariant vector bundle can be equivariantly quantized.}

\section{Further Remarks on Bundles}
\subsection{Uniqueness}
\label{unique}
Equivariant bundles and modules are classified by equivariant 
$K$-theory. The equivariant vector bundles over $\Orl$ all have 
equivalence classes in $K^0_G(\Orl)$. As has been mentioned 
(Sec.~\ref{equi_bun}), these bundles are classified by 
representations of $H$. From this it is easy to see that an 
equivariant bundle is uniquely specified by its $K$-class. Similarly, 
an equivariant module of $\AA$ is uniquely specified by its $K$-class 
in $K^G_0(\AA)$. The equivariant general quantizations of vector 
bundles are equivariant modules of $\AA$, and are thus classified by 
$K^G_0(\AA)$. Since $\C(\Orl)=\AA/\AA_0$, there is a corresponding 
six-term periodic exact sequence in $K$-theory. Part of this sequence 
reads
\beq
K^G_0(\AA_0)\to K^G_0(\AA) \to K_G^0(\Orl) \to K^G_1(\AA_0)
\;.\label{Kexact}\eeq

The ideal $\AA_0=\Gn(\A_\N)$ is the \cs-direct sum of the algebras 
$\A_N$; therefore $K^G_*(\AA_0)=\bigoplus_{N=1}^\infty K^G_*(\A_N)$. 
Because $\A_N$ is the matrix algebra on a simple representation of 
$G$, its equivariant $K$-theory is very simple. In degree $0$, 
$K^G_0(\A_N)= {\mathcal R}(G)$ the unitary representation ring of 
$G$. In degree $1$, $K^G_1(\A_N)=0$. This simplifies the exact 
sequence \eqref{Kexact}. Now it reads
\beq
{\mathcal R}(G)^{\oplus\infty} \to K^G_0(\AA) \onto K_G^0(\Orl) \to 0 
\;.\eeq

Firstly, this shows that --- at the level of $K$-theory --- any 
equivariant bundle has an equivariant quantization, since it has a 
preimage in $K^G_0(\AA)$. This corroborates the conclusion of Section 
\ref{identification}. Secondly, this describes the variety of 
possible quantizations of a given bundle. If two equivariant 
$\AA$-modules quantize the same bundle, then the difference of their 
$K$-classes is in the image of ${\mathcal R}(G)^{\oplus\infty}$, but 
that is an algebraic direct sum; it consists of sequences with only
finitely many nonzero terms, and each term concerns a single $N$. 
This means that if both $V_\N$ and $V'_\N$ are quantizations of $V\!$, 
then for all $N$ sufficiently large, $V_N\cong V'_N$.

Given this conclusion, the choice of $V_N$'s in Eq.~\eqref{Vl_form} must be 
the unique one given by a simple formula.

\subsection{Geometric Quantization}
\label{geoQ}
For each $N$, the fundamental module $(N\Lambda)$ of $\A_N$ is of 
course a module. It is tempting to ask if these together form the 
quantization of some bundle, but they do not. The $\AA$-module formed 
by assembling these is not projective.

It is reasonable to instead ask --- separately for each $N$ --- what 
bundle's equivariant quantization (by the construction of 
Sec.~\ref{Q_vec_buns}) has $V_N=(N\Lambda)$? This is easily answered;
\beq
(N\Lambda)=(N\Lambda)\otimes(0)=(N\Lambda)\otimes(N\Lambda^*-N\Lambda^*) 
=V^{-N\Lambda^*}_N
\;.\eeq
Using the identity that $[N\Lambda]^*=[-N\Lambda]$ (see 
App.~\ref{coadjoint_classification}), the corresponding bundle is 
$V^{-N\Lambda^*}=[N\Lambda]\times_H G$. 

The $H$-representation $[N\Lambda]$ is one-dimensional; this bundle 
is therefore of rank 1 --- i.~e., it is a line bundle. In geometric 
quantization of $\Orl$, the fundamental module $(N\Lambda)$ of $\A_N$ 
is constructed as the space of holomorphic sections of this very line 
bundle. 

\subsection{Bimodules}
\label{bimodules}
For a commutative algebra, any module can automatically be considered 
a bimodule; simply define right multiplication to be equal to left 
multiplication. However, it is not generally the case that when a 
vector bundle is quantized, the corresponding module continues to be 
a bimodule. The right side (the row-vector factor) of $V^\lambda_N$ 
is $(N\Lambda^*+\lambda)$ and does not in general admit any 
equivariant right multiplication by $\A_N\equiv\End(N\Lambda)$.

If $V_N$ is an $\A_N$-bimodule, then it must contain a factor of 
$(N\Lambda)$ to accommodate the left multiplication, and a separate 
factor of $(N\Lambda^*)$ to accommodate the right multiplication. It 
must therefore be the tensor product of $\A_N$ itself by some 
representation. The corresponding classical bundle is then the 
trivial bundle with fiber equal to that representation. This is an 
unpleasantly restrictive class. 

A slightly broader class of bundles results if we allow the quantum 
modules to be multiplied from the left and right by \emph{different} 
$\A_N$'s. This is enough to make $\VV$ an $\AA$-bimodule, and is also 
contrary to the philosophy of each $N$ being a separate step along 
the way to the classical limit. The irreducibles of this class of 
modules are of the form 
\beq
V_N=(N\Lambda)\otimes([N+m]\Lambda)\otimes(\lambda) = 
V^{m\Lambda^*}_N\otimes(\lambda)
\;.\eeq 
The corresponding classical bundles are a slightly more interesting 
class than trivial bundles, but still quite restrictive. This can be 
extended a little further in some cases by using a larger parameter 
set $\I$. It remains to be seen whether this class of modules is 
useful. 

\section{The Case of the 2-sphere}
\label{twosphere}
The group $\sut$ is the most elementary compact, simple Lie group, so 
the simplest example of what has been described here is for $G=\sut$. 
There is only one distinct coadjoint orbit for $\sut$; it is the 
2-sphere. As a coset space $S^2=\sut/\Uo$.

The positive Weyl chamber of $\sut$ is $\C_+=\Rplus$. Thought of as the 
parameter space for $S^2$'s, this is the set of radii. In deference 
to standard physics notation, I will identify the dominant weights 
with positive half-integers. The irreducible representations are thus 
$(0)$, $(\frac12)$, $(1)$, \emph{et cetera}. The most appropriate 
choice for $\Lambda$ is $\frac12$.

The Lie algebra $\mathfrak{su}(2)$ is generated by $J_1$, $J_2$, and 
$J_3$, with the commutation relations $[J_i,J_j]_- = 
i\epsilon^k_{\;ij} J_k$; that is, $[J_1,J_2]_-=iJ_3$, \emph{et 
cetera}.

The standard choice for the Cartan subalgebra $\Ca$ is the 
one-dimensional span of $J_3$. The weights are just the eigenvalues 
of $J_3$. In the representation $(\frac{N}2)$ the highest weight 
vector satisfies $J_3\Psi^{N/2}=\frac{N}2\Psi^{N/2}$.

There is a single (quadratic) Casimir operator $\C_1(J) = J^2 \equiv 
J_1^2+J_2^2+J_3^2$. Its eigenvalue on the representation 
$(\frac{N}2)$ is $\frac{N}2[\frac{N}2+1]$. 

There is a single Serre relation for $\End(\frac{N}2)$. In terms of 
the element $J_+ := \tfrac12 (J_1+iJ_2)$, the relation is that 
$J_+^{N+1}=0$. Although this is expressed in a noninvariant way, this 
condition really is invariant; it could equivalently be expressed in 
terms of many other possible combinations of 
$J$'s. 

The logic of the Serre relation is that the representation 
$(\frac{N}2)$ is $N+1$ dimensional. It can be decomposed into 
one-dimensional weight subspaces ($J_3$ eigenspaces). The operator 
$J_+$ shifts these weight subspaces; it maps the subspace with weight 
$m$ to the subspace with weight $m+1$ (the next higher possible 
weight). $J_+$ can be applied to some $J_3$-eigenvector no more than 
$N$ times before there are no more eigenvalues available, and the 
result must be $0$. Therefore $J_+^{N+1}$ applied to anything in 
$(\frac{N}2)$ must give $0$.
\gap

We can construct a general quantization by the method of Section 
\ref{gen_pic}. The generators $x_i:=\Po(N^{-1} J_i)$ of the resulting 
$\Ai$ satisfy the relations of commutativity and 
$x_1^2+x_2^2+x_3^2=\frac14$. Obviously this shows $\Ai$ to be the 
continuous functions on the sphere of radius $\frac12$ in 
$\mathfrak{su}(2)^*\cong\R^3$.

All $\sut$-representations are self-dual. Because of this, the 
constructions of $i_N$ and $p_N$ are even closer than in the general 
case in Section \ref{coherent}. Decompose the tensor product 
$(\frac12)\otimes(\frac{N}2)=(\frac{N+1}2)\oplus(\frac{N-1}2)$. There 
is a representation of $\A_N$ on this that acts trivially on the 
$(\frac12)$ factor; for an element $a\in\A_N$, the $(\frac{N+1}2)$ 
corner of this representation matrix is $i_N(a)\in\A_{N+1}$; the 
$(\frac{N-1}2)$ corner is $p_N(a)\in\A_{N-1}$.

In this case it is possible to construct a simple and (partly) 
explicit formula for the ``product'' kernel $K_N$. The key is to use 
the identification $S^2=\co \mathrm{P}^1 = \PP(\frac12)$. The 
geodesic distances on $S^2$ are given by the Fubini-Study metric; for 
two points $[\psi],[\varphi]\in\PP(\frac12)$ the distance 
$d_{S^2}([\psi],[\varphi])$ is determined by 
$$
\cos^2 \left[d_{S^2}\left([\psi],[\varphi]\right)\right] := 
\frac{\;\left|\left<\psi|\varphi\right>\right|^2}%
{\left<\psi|\psi\right>\left<\varphi|\varphi\right>}
\;.$$
Recall that this is the sphere of radius $\frac12$, so $0\leq 
d_{S^2}(x,y)\leq\frac\pi2$. Comparing this with the formulas 
\eqref{kernel} or \eqref{K_unnorm} for $K_1$ gives that
$$
\lvert K_1(x,y,z)\rvert = 4\cdot \cos [d_{S^2}(x,y)] \cos 
[d_{S^2}(y,z)] \cos [d_{S^2}(z,x)] 
\;.$$ 
Noting \eqref{convenient}, this gives for arbitrary $N$ that
$$
\lvert K_N(x,y,z)\rvert = (N+1)^2 \cos^N [d_{S^2}(x,y)] \cos^N 
[d_{S^2}(y,z)] \cos^N [d_{S^2}(z,x)] 
\;.$$ 
What remains to be determined is the phase. This has no simple 
formula, but is easily understood geometrically: $\arg K_N(x,y,z)$ is 
$2N$ times the area of the geodesic triangle on $S^2$ with vertices 
$x$, $y$, and $z$. To see this, show that this is true for 
infinitesimal triangles and that this quantity is additive when a 
triangle is decomposed into smaller triangles. Clearly $K_N$ does 
indeed  become sharply peaked as $N\to\infty$.
\gap

Since the isotropy group of $S^2$ is $H=\Uo$, the classification of 
equivariant vector bundles over $S^2$ is extremely simple. The 
irreducible bundles are classified by irreducible representations of 
$\Uo$, which are in turn indexed by half integers. Denote these 
representations by $[m]$ for any $m\in\frac12\Z$. Since these 
representations are all one-dimensional, the irreducible bundles are 
all rank-one. 

Under the restriction $\sut\hookleftarrow\Uo$, an irreducible 
representation of $\sut$ decomposes into a direct sum of irreducible 
$\Uo$-representations. This is simply 
$$
(j)\to[-j]\oplus[-j+1]\oplus\dots\oplus[j]
$$
for any $j\in\frac12\Z$.

Let $W^m:=[m]\times_\Uo \sut$ be the equivariant vector bundle over 
$S^2$ with fiber $W^m_o=[m]$. The space of continuous sections 
$\G(W^m)$ is a completion of the space of polynomial sections 
$\Gp(W^m)$. As an $\sut$-representation, $\Gp(W^m)$ is a direct sum 
of irreducibles and is easily computed. The representation $(j)$ 
occurs precisely once in $\Gp(W^m)$ if and only if $[m]$ occurs in 
the decomposition of $(j)$; in other words, when $j\equiv m \mod 1$ 
and $j\geq\lvert m\rvert$. So,
$$
\Gp(W^m) = (\lvert m\rvert)
\oplus(\lvert m\rvert+1)\oplus(\lvert m\rvert+2)\oplus\cdots
\;.$$

So, what is the quantization of $W^m$? We need to know for which $j$ 
does $W^m=V^j$. Equation \eqref{the_point} shows that $[m]=[j^*]^*$. 
Since all $\sut$-representations are self-dual, $j^*=j$. All 
irreducible $\Uo$-representations are one-dimensional, so (see 
App.~\ref{coadjoint_classification})  $[j]^*=[-j]$. This means that 
$j=-m$, and the quantization of $W^m$ is $\VV^{-m}$.

As an $\sut$ representation 
$$
V_N^{-m}=(\tfrac{N}2)\otimes(\tfrac{N}2-m) = 
(\lvert m\rvert)\oplus(\lvert m\rvert +1)\oplus\dots\oplus(N-m)
\;.$$ 
Clearly, modulo completion, $V_N^{-m}$ in the limit $N\to\infty$ 
becomes the same $\sut$-representation as $\G(W^m)$. This applies in 
particular when $\lambda=0$ so $V_N^0=\A_N$ and $\G(W^0)=\C(S^2)$. 
With a little more work, this consistency check can be carried out 
for all the complex projective spaces $\co \mathrm{P}^n$.
\gap

In a way, it may seem odd to be using $\sut$ as the symmetry group 
for $S^2$. The group of \emph{distinct} orientation preserving 
isometries of $S^2$ is $\SO(3)$; the group $\sut$ is its 
simply-connected, double cover. If we had used the smaller group, we 
would have artificially excluded all the $\A_N$'s with $N$ odd. 
Although $\SO(3)$ acts on all the algebras $\A_N$, we need the 
$\sut$-representation $(N\Lambda)$ in order to construct $\A_N$. 
Another reason is that many of the vector bundles on $S^2$ are 
$\sut$-equivariant, but not $\SO(3)$-equivariant.

It is generally the case that the simply connected $G$ is not the 
minimal symmetry group of a coadjoint orbit. Indeed, the minimal 
symmetry group of $\Orl$ is the group $G'=G/Z(G)$ (the ``adjoint 
group'') which maximizes the fundamental group $\pi_1(G')$. 
Nevertheless, the simply connected $G$ is the easiest to deal with, 
and most fruitful, choice. 

\section{Final remarks}
\label{conclusions}
One motivation for considering the limit quantization approach for 
bundles comes from physics. If this sort of quantization is used as a 
regularization technique, then it would be desirable to do a 
``renormalization group'' analysis. This involves going from one 
level of regularization to a coarser one with fewer degrees of 
freedom. In order to do this we need a sort of coarse-graining map 
that associates a given field configuration with a coarser field 
configuration, ignoring some of the degrees of freedom of the 
original. 

In $n$-dimensional lattice regularization, the space is approximated 
by a lattice. The coarse-graining is accomplished by grouping the 
lattice points into groups of $2^n$ and averaging the field values at 
those $2^n$ points. This field value is then given to a single point 
of the new, coarsened lattice which has $2^{-n}$ times as many 
points. The degrees of freedom are thus reduced (drastically) by a 
factor of $2^n$.

Classically, field configurations are sections of vector bundles. If 
quantization is used as a regularization technique, the field 
configurations are the vectors in the quantum modules $V_N$. 
Coarse-graining means going from $N$ to $N-1$. The coarse-graining 
map is $\pi_N$. The degrees of freedom vary as only a polynomial 
function of $N$, so $\dim V_N/\dim V_{N-1} \approx 1$ for large $N$. 
This is far gentler than lattice regularization. I hope to discuss 
this, and related matters in a future paper.
\gap

Another reason for using constructions in terms of limits, as I have 
here, is simply that it is the most convenient approach when dealing 
with coadjoint orbits. When dealing with the quantization of a more 
general symplectic manifold, objects such as the Hilbert space $\Hi_N$ 
are constructed as spaces of sections over the manifold; everything is 
constructed from the manifold. In the case of coadjoint orbits, 
however, $\Hi_N$ is constructed directly as a $G$-representation. We can 
actually deal more explicitly with the algebra $\A_N$ than with the 
algebra $\C(\Orl)$. For this reason, it is more convenient to 
construct the classical structures from the quantum structures, rather 
than \emph{vice versa}.
\gap

The construction of the maps $I_N$ and $P_N$ in Section \ref{coherent} 
is standard \cite{ber,g-p,per}. In the terminology of Berezin 
\cite{ber}, $P_N(a)$ is the contravariant symbol of $a$, and an 
element of the preimage $I_N^{-1}(a)$ is a covariant symbol of $a$.
\gap

The idea of directed limit quantization here is based on a 
construction by Grosse, Klim\v{c}\'{\i}k, and Pre\v{s}najder in 
\cite{g-k-p2}. In that case the quantization of the $S^2$ was being 
discussed. Their choice of $i_N$ is different and is based on the 
criterion of preserving the $L^2$-norm from one algebra to the next. 
My choice as based on the criterion of compatibility with the 
standard $I_N$'s. It can easily be checked that $I_N$ never preserves 
the $L^2$ norms, and therefore my choice of $i_N$'s never satisfies 
their  criterion. 
\gap

In \cite{g-k-p3}, Grosse, Klim\v{c}\'{\i}k, and Pre\v{s}najder 
constructed quantized vector bundles for the special case of $S^2$. 
Their result is the same as mine for that case (see 
Sec.~\ref{twosphere}). 
\gap

To reiterate, the main conclusion of this paper is that when the 
coadjoint orbit $\Orl=G/H$ through $\Lambda$ is quantized to give a 
sequence of matrix algebras
$$
\A_N = \End(N\Lambda)
\;,$$
the equivariant vector bundle
$$
V^\lambda = [\lambda^*]^* \times_H G
\eqno\eqref{the_point}$$
quantizes to a corresponding sequence of $\A_N$-modules
$$
V^\lambda_N = (N\Lambda)\otimes(N\Lambda^*+\lambda)
\;.\eqno\eqref{Vl_form}$$
\gap

In \cite{haw4} I will continue by describing analogous results 
in the more general case of compact K\"ahler manifolds.

\begin{appendix}
\section{Sections}
\label{sections}
Before discussing the construction of limits, it is worthwhile to 
clarify the notations for different spaces of sections of the bundles of 
algebras and modules.
Given a noncompact base space, there are several useful types of 
continuous sections of a vector bundle, all of which are equivalent 
for a compact base space. For the base space $\N$, sections are the 
same thing as sequences. For legibility, I will sometimes write 
sections as sequences in that case.

The space of all continuous sections of a vector bundle $E$ is 
denoted $\G(E)$. If $E$ is a bundle of algebras, then $\G(E)$ is an 
algebra. However, for a bundle $\A_\I$ of \cs-algebras $\G(\A_\I)$ is 
not a \cs-algebra since the sup-norm diverges. For a discrete base space, 
this is the algebraic direct product.

The space of continuous sections with compact support is denoted 
$\Gc(E)$. For a bundle of algebras, this is an ideal inside $\G(E)$. 
For the \cs-bundle $\A_\I$ this space $\Gc(\A_\I)$ has a \cs-norm, 
but is not complete and therefore not \cs. For a discrete base space, 
this is the algebraic direct sum.

If the fibers of $E$ are normed (as \cs-algebras are), then two more 
types of section can be defined. $\Gb(E)$ is the space of sections of 
bounded norm. For the \cs-bundle $\A_\I$, $\Gb(\A_\I)$ is a 
\cs-algebra; the norm of a section is the supremum of the norms at 
all points of $\I$. For \cs-algebras over a discrete base space this 
is the \cs-direct sum.

$\Gn(E)$ is the space of sections such that the norms converge to $0$ 
approaching $\infty$. To be precise, any arbitrarily low bound on the 
norms is satisfied on the complement of some compact set. This is the 
norm closure of $\Gc(E)$. For the \cs-bundle $\A_\I$, $\Gn(\A_\I)$ is 
a closed ideal in $\Gb(\A_\I)$. For \cs-algebras over a discrete base 
space, this is the \cs-direct product.

These spaces of sections are related by 
$\Gc\subset\Gn\subset\Gb\subset\G$.

The appropriate notion of a bundle of \cs-algebras is that of a 
continuous field of \cs-algebras. This is discussed extensively in 
\cite{dix}.

\section{Limits}
\label{lims}
\subsection{Direct limit of algebras}
\label{alg_dir_lim}
Since we are assuming the index set to be $\N$, sections of $\A_\N$ 
can also be thought of as sequences. 

In the category of vector spaces, the limit of a directed system of 
algebras is 
\beq
\veclim\{\A_*,i_*\} := \oAA/\Gc(\A_\N)
\label{veclimdef}\;,\eeq
where
\beq
\oAA:=\left\{a\in\G(\A_\N) \bigm\vert \exists M\;\forall N\geq M: 
a_{N+1} = i_N(a_N)\right\}
\label{oAA_def}\;.\eeq
The injections $i_N$ are meant to identify $a_N$ to $i_N(a_N)$; 
\eqref{veclimdef} therefore gives the set of sequences which for 
sufficiently large $N$ become constant, modulo the sequences which 
for sufficiently large $N$ are $0$. Thinking of 
$\A_N\subset\A_{N+1}$, the limit is heuristically the union 
$\bigcup_{N\in\N}\A_N$ of this nested sequence.

Usually, one works in the category of \cs-algebras in which the 
morphisms are $*$-homomorphisms. If the $i_N$'s are 
assumed to be $*$-homomorphisms, then the \cs-algebraic limit (see 
\cite{fil}) of finite-dimensional algebras will be (by definition) an 
AF-algebra. This is far too restrictive a class of algebras in this 
context; a commutative AF-algebra is isomorphic to the continuous functions 
on a totally disconnected, zero-dimensional space (see \cite{weg}). In 
order to avoid this restriction, we must allow the $i_N$'s to be some 
more general type of morphism. Firstly, these must be linear, and I 
will assume (perhaps unnecessarily) that they are unital (i.~e., 
$i_N(1)=1$). 

Several convergence conditions on the $i_N$'s will also be needed. 
The first condition is that the $i_N$'s be norm-contracting maps; 
this means $\forall a\in\A_N$, $\norm{i_N(a)}\leq\norm{a}$. There is 
a fairly nice class of norm-contracting maps for \cs-algebras; these 
are the completely positive maps (see \cite{lan}). All of the $i_N$'s 
and $p_N$'s constructed in this paper are completely positive; 
however, I am not relying on that property in general. The 
norm-contracting condition ensures $\oAA\subset\Gb(\A_\N)$. 

Since each $\A_N$ is a \cs-algebra, each has a \cs-norm. The natural 
norm on the limit is the limit of these; that is, for any equivalence class 
$[a]\in\veclim\{\A_*,i_*\}$ define 
\beq
\norm{[a]} := \lim_{N\to\infty} \norm{a_N}
\label{norm_form}\;.\eeq 
The norm-contracting condition guarantees that this is well defined, 
since it is a limit of a sequence that is (for sufficiently large 
$N$) strictly nonincreasing and bounded from below (by $0$). 

To ensure that this is truly a norm requires a second condition --- 
that it be nondegenerate. That is, $a\neq0 \implies \norm{a}\neq0$. 
This is equivalent to the condition that $\oAA\cap\AA_0=\Gc(\AA_\N)$, 
where $\AA_0=\Gn(\A_\N)$.

This means that $\Gc(\A_\N)$ can be replaced by $\AA_0$ in 
\eqref{veclimdef}, and $\oAA$ is naturally embedded in the 
\cs-algebra $\Gb(\A_\N)/\AA_0$. The norm \eqref{norm_form} agrees 
with the natural norm on this quotient. Now define 
$\Ai=\dlim\{\A_*,i_*\}$ as the closure of $\oAA/\AA_0$ in 
$\Gb(\A_\N)/\AA_0$, or equivalently as the abstract norm completion 
of $\oAA/\AA_0$.

Also define $\AA\subset\Gb(\A_\N)$ as the norm closure of 
$\oAA\subset\Gb(\A_\N)$. Another construction of $\Ai$ is 
$\Ai=\AA/\AA_0$; this shows that if we view sections in $\Gb(\A_\N)$ 
as sequences, then $\AA$ is the subspace of sequences which converge 
into $\Ai$.

It is not \emph{a priori} true that $\Ai$ is an algebra; this 
requires a third (and final) condition. Require that $\Ai$ be 
algebraically closed in $\Gb(\A_\N)/\AA_0$. This is equivalent to 
requiring that $\AA\subset\Gb(\Ai)$ be algebraically closed.

Assuming these conditions, both $\Ai$ and $\AA$ are norm closed 
subalgebras of \cs-algebras; they are therefore \cs-algebras 
themselves.

For each $N$, there is a canonical injection 
\beq
I_N:\A_N\into\veclim\{\A_*,i_*\}\subset\Ai
\eeq
which takes $a\mapsto \left[\left(0,\ldots,0,a,i_N(a),i_{N+1}\circ 
i_N(a),\ldots\right)\right]$. Heuristically, $I_N = \ldots\circ 
i_{N+1}\circ i_N$.
\gap

If we are trying to prove that a given directed system $\{\A_*,i_*\}$ 
truly converges to a given $\Ai$, the third convergence condition is 
the most critical. Using the notation 
$i_{N,M}:=i_{M-1}\circ\ldots\circ i_N:\A_N\into\A_M$, an equivalent 
statement is that $\forall N \; \forall a,b\in\A_N$,
\beq
\lim_{m\to\infty} I_{N+m}[i_{N,N+m}(a)\: i_{N,N+m}(b)] = I_N(a) I_N(b)
\;.\eeq

If there are left inverses $I^{\text{inv}}_N$ (such that 
$I^{\text{inv}}_N\circ I_N = \id$), chosen so that the sections 
$N\mapsto I^{\text{inv}}_N(f)$ are continuous, then there is a 
simpler statement. This convergence condition becomes $\forall 
f_1,f_2\in\Ai$ 
\beq
 I_N[I^{\text{inv}}_N(f_1)I^{\text{inv}}_N(f_2)] 
\toas f_1 f_2
\;.\eeq 
This is the form used in Section \ref{convergence}.

In this circumstance it is also necessary to check that 
$\veclim\{\A_*,i_*\}\subset\Ai$ really is dense. This means that 
$I_N$ needs to be ``asymptotically onto''\negthinspace. Using the 
$I^{\text{inv}}_N$'s, this simplifies to the requirement that 
$\forall f\in\Ai$, $I_N\circ I^{\text{inv}}_N(f) \to f$ as 
$N\to\infty$.

Although this was done for the index set $\N$, it can trivially be 
generalized to any directed set.

\subsection{Direct limit of modules}
\label{mod_dir_lim}
Given a directed system $\{V_*,\iota_*\}$ of finitely generated, 
projective modules of each $\A_N$, we would like to construct a limit 
module of the limit algebra $\Ai$. The construction must work in the 
special case that the system is just $\{\A_*,i_*\}$. The vector space 
direct limit $\veclim\{V_*,\iota_*\}$ is not itself an $\Ai$-module; 
it needs to be completed somehow. Completion is usually done with 
some norm, but there is generally no natural norm on the $V_N$'s. 
Instead, complete algebraically.

The algebraic direct product $\G(V_\N)$ is a $\G(\A_\N)$-module, and 
by restriction an $\AA$-module. From the construction of the vector 
space direct limit, start with the vector space 
\beq
\oVV := \left\{\psi\in\G(V_\N) \bigm\vert \exists M\;\forall N\geq M: 
\psi_{N+1} = \iota_N(\psi_N)\right\}
\;.\eeq
Now define $\VV$ as the span of $\AA\oVV\subset\G(V_\N)$. I insist 
that $\VV$ be a finitely generated $\AA$-module, so there is a 
convergence condition that any element of $\VV$ can be written as the 
sum of a bounded number of elements of $\AA\oVV$. In other words, 
$\AA\oVV+\dots+\AA\oVV$ stabilizes for some finite number of summands.

It is now easy to construct an $\Ai$-module. The ideal $\AA_0$ 
induces a submodule $\AA_0\!\VV\subset\VV\!$, and the quotient 
$V_\infty:=\VV/\AA_0\!\VV$ is an $\Ai$-module. This is the direct 
limit of modules. Note that its construction requires the map 
$\Po:\AA\onto\Ai$ but does not require any other quantization 
structure for the algebras.

\subsection{Inverse limit of algebras}
\label{alg_inv_lim}
The limit of the inverse system of algebras is easier to construct. 
It is
\beq
\Ai=\ilim\{\A_*,p_*\} := \left\{a\in\Gb(\A_\N) \bigm\vert \forall N: 
a_{N-1}=p_N(a_N)\right\}
\;.\eeq

Again, the $p_N$'s should not be required to be homomorphisms, and 
again, convergence conditions are necessary.

This limit also inherits a norm $\norm{a}:=\lim_{N\to\infty} 
\norm{a_N}$. This is well defined if the $p_N$'s are required to be 
norm-contracting. It is then the limit of a nondecreasing sequence 
that is bounded from above. No additional condition is required to 
make this nondegenerate since $\norm{a}\geq\norm{a_N}$.
$\Ai$ is already complete with respect to this norm.

Since $\Ai$ consists of sequences of nondecreasing norm, the 
intersection with $\AA_0=\Gn(\A_\N)$ is $0$. This means that $\Ai$ 
injects naturally into $\Gb(\A_\N)/\AA_0$. Define $\AA$ to be the 
preimage of $\Ai$ by the quotient homomorphism 
$\Gb(\A_\N)\onto\Gb(\A_\N)/\AA_0$; this gives $\AA=\Ai + 
\AA_0\subset\Gb(\A_\N)$.

This $\Ai$ is also not \emph{a priori} an algebra. We again need the 
condition that $\Ai\subset\Gb(\A_\N)/\AA_0$ be algebraically closed. 
This is equivalent to requiring that $\AA\subset\Gb(\A_\N)$ be 
algebraically closed. If $\Ai$ and $\AA$ are algebraically closed, 
then they are \cs-algebras.

For each $N$, there is a canonical surjection
$P_N:\ilim\{\A_*,p_*\} \onto \A_N$ which simply takes $a\mapsto a_N$. 
Heuristically $P_N = p_{N+1}\circ p_{N+2}\circ\ldots$.
\gap

This last convergence condition is again the most critical. If we are 
testing whether $\ilim\{\A_*,p_*\}\stackrel?= \Ai$, then an 
equivalent statement is $\forall f_1,f_2\in\Ai$,
\beq
\lim_{N\to\infty} \norm{P_N(f_1)P_N(f_2) - P_N(f_1 f_2)} = 0
\;.\eeq

If there are right inverses $P^{\text{inv}}_N$ (such that $P_N\circ 
P^{\text{inv}}_N = \id$), chosen so that $P^{\text{inv}}_N\circ 
P_N(f)\to f$ as $N\to\infty$, then there is a simpler statement. This 
convergence condition becomes $\forall f_1,f_2\in\Ai$,
\beq
P^{\text{inv}}_N[P_N(f_1)P_N(f_2)] \toas f_1 f_2
\;.\eeq 
This is the form used in Section \ref{convergence}.

\subsection{Inverse limit of modules}
\label{mod_inv_lim}
This construction is very much the same as in \ref{mod_dir_lim} for a 
direct limit of modules. For an inverse system $\{V_*,\pi_*\}$ of 
modules, first construct the vector space 
\beq
\oVV := \left\{\psi\in\G(V_\N) \bigm| \forall N : 
\psi_{N-1}=p_N(\psi_N)\right\}
\;.\eeq
Again define $\VV$ as the span of $\AA\oVV\!$, and the convergence 
condition is that $\AA\oVV+\dots+\AA\oVV$ stabilizes for some 
finite number of summands.

Define $\ilim\{V_*,\pi_*\}:=\AA/\AA_0\!\VV$.

\section{Review of Representation Theory}
\label{representation_theory}
Let $G$ be a compact, simply connected, semisimple Lie group. This 
always contains a Cartan subgroup $T$. This is a maximal abelian 
subgroup which is always of the form $\Uo^\ell$ (a torus group). Any 
two Cartan subgroups of $G$ are conjugate, so it is irrelevant which 
one we now fix and call \emph{the} Cartan subgroup. Since the 
irreducible representations of $\Uo$ are one-dimensional and 
classified by $\Z$, the irreducible representations of $T$ are 
one-dimensional and classified by the lattice $\Z^\ell$. 

The \emph{Cartan subalgebra} $\Ca\subset\g$ is the Lie algebra of 
$T$. Any vector in an irreducible representation of $T$ is an 
eigenvector of any element of $\Ca$; the eigenvalue depends linearly 
on the position of the representation in the above lattice (and on 
the element of $\Ca$). The lattice is therefore naturally thought of 
as lying in the dual $\Ca^*$ of the Cartan subalgebra. It is called 
the \emph{weight lattice}. There is a natural inner product (the 
Cartan-Killing form) on the 
Lie algebra $\g$; using this, there is a natural sense in which 
$\Ca^*\subset\g^*$. 

There are some symmetries to $\Ca^*$, residual from the action of $G$ 
on $\g^*$. The symmetry group of $\Ca^*$ is the subgroup of $G$ that 
preserves $\Ca^*\subset\g^*$, modulo the subgroup that acts trivially 
on $\Ca^*$. This is called the \emph{Weyl group} $W$ and is finite. 
Since both are naturally constructed from the pair $\Ca\subset\g$, 
the Weyl group preserves the weight lattice. The Weyl group is 
generated by a set of reflections across hyperplanes in $\Ca^*$. 
These plains divide $\Ca^*$ into wedges called \emph{Weyl chambers}; 
each Weyl chamber is a fundamental domain of the $W$ action on 
$\Ca^*$, this means that the $W$-orbit of any point of $\Ca^*$ 
intersects a given closed Weyl chamber at least once and intersects 
the interior of a given Weyl chamber at most once. 

We can choose a basis of the weight lattice; that is, a set of 
\emph{fundamental weights} $\{\pi_j\}$ such that the weight lattice 
is the integer span $\sum_j \Z\pi_j$. Given the choice of $\Ca$, the 
fundamental weights are unique modulo the freedom to change their 
signs. Fix a set of fundamental weights. The natural index set for 
the fundamental weights is the set of vertices of the Dynkin diagram of 
$\g$. The positive span of the fundamental weights $\sum_j \Rplus\pi_j$ 
is precisely a (closed)  Weyl chamber. Call this \emph{the positive 
Weyl chamber} $\C_+$. The weights that lie in $\C_+$ are nonnegative 
integer combinations of the fundamental weights and are called 
\emph{dominant} weights. Since it is a fundamental domain of the $W$ 
action, the positive Weyl chamber $\C_+$ can naturally be identified 
with $\Ca^*/W$.
\gap

Given an irreducible representation of $G$, it can also be regarded as a 
$T$ representation. The representation space therefore naturally 
decomposes into a direct sum of subspaces associated with different 
weights. The set of weights that occur is $W$-invariant. The subspace 
associated with the dominant weight furthest from $0$ is always 
1-dimensional; that weight is called the \emph{highest weight} of the 
representation. Nonisomorphic irreducible representations have 
distinct highest weights and any dominant weight is the highest 
weight of some representation. The irreducible representations of $G$ 
are therefore exactly classified by dominant weights. I denote the 
representation space with highest weight $\lambda$ as $(\lambda)$.

Weights are additive under the tensor product. If two 
vectors have weights $\lambda$ and $\mu$, then their tensor product 
has weight $\lambda+\mu$. Because of this, the highest weight of the 
(reducible) representation $(\lambda)\otimes(\mu)$ is $\lambda+\mu$. 
The decomposition of $(\lambda)\otimes(\mu)$ into irreducibles will 
therefore always contain precisely one copy of $(\lambda+\mu)$; this 
irreducible representation is called the Cartan product of 
$(\lambda)$ and $(\mu)$.

For each irreducible representation, we can choose a normalized 
vector $\Psi^\lambda\in(\lambda)$ in the highest weight subspace. 
This is called a \emph{highest weight vector}. Their phases are 
arbitrary, but can be chosen consistently so that 
$\Psi^\lambda\otimes\Psi^\mu=\Psi^{\lambda+\mu}\in(\lambda+\mu)\subset 
(\lambda)\otimes(\mu)$. 
\gap

The linear dual of an irreducible representation is also an 
irreducible representation; we can therefore define $\lambda^*$ by 
the property $(\lambda^*)=(\lambda)^*$. This is a linear 
transformation on the weight lattice; it simply permutes the 
fundamental weights and is given by an automorphism (possibly 
trivial) of the Dynkin diagram. 

Whenever $\lambda-\mu$ is a dominant weight, 
$(\lambda)\otimes(\mu^*)$ will contain precisely one copy of 
$(\lambda-\mu)$. In particular, if $\lambda=\mu$ this says that 
$(\lambda)\otimes(\lambda^*)$ contains one copy of the trivial 
representation; this is little more than the definition of the dual.

\section{Coadjoint Orbits}
\label{restrictions}
The purpose of this appendix is to describe the rationale for 
restricting attention to coadjoint orbits, and then to discuss some 
of the structure of coadjoint orbits. Toward this goal, I first 
discuss a more general structure:

\subsection{Symplectic structure}
\label{symplectic}
Thus far I have entirely avoided mentioning something which is 
usually mentioned first in discussions of quantization --- the 
symplectic structure. 

Assume $\M$ to be a manifold. Suppose that part of our quantization 
structure is a system of maps $I_N:\A_N\into\C(\M)$, identifying 
quantum operators to classical functions. We can choose a system of 
right inverses; that is, maps $I^{\text{inv}}_N:\C(\M)\onto\A_N$ such 
that $I_N\circ I^{\text{inv}}_N$ is the identity map $\A_N\to\A_N$. 
Using these, we can pull back the products on each of the $\A_N$'s to 
$\C(\M)$, giving a sequence of products 
\beq
f*^N f' = I_N[I^{\text{inv}}_N(f) I^{\text{inv}}_N(f')]
\label{star}\;.\eeq 
By construction, these converge to the ordinary product of functions 
as $N\to\infty$. 

Suppose that the quantization is compatible with the smooth structure 
of $\M$ in the sense that for smooth functions $f,f'\in\C^\infty(\M)$ 
the correction $f*^N f' - ff'$ is of order\footnote{This can be 
generalized slightly by replacing $\frac1N$ with some other function 
$\hbar(N)$ that goes to $0$ as $N\to\infty$. The implication 
(existence of Poisson bracket) remains the same.} $\frac1N$. I will 
assume that any quantization of interest satisfies this.
This compatibility means that the function
\beq
\{f,f'\}:= \lim_{N\to\infty} -iN\left(f*^N f' - f'*^N f\right)
\eeq
is well defined. This is the Poisson bracket of $f$ and $f'$; it is 
easily seen to be, by construction, antisymmetric and a derivation in 
both arguments. This means that there exists an antisymmetric, 
contravariant, rank-2 tensor\footnote{This is also called a 
bivector.} $\pi^{ij}$ such that the Poisson bracket is given by 
$\{f,f'\} = \langle\pi,df\wedge df'\rangle \equiv \pi^{ij} df_i\, 
df'_j $. 

With the assumption that the algebras $\A_N$ are finite-dimensional, 
the $\pi$ should be nondegenerate if thought of as a map from 1-forms 
to tangent vectors. This means that it has an inverse 
$\omega=\pi^{-1}$, which is naturally a 2-form.
The Poisson bracket also satisfies the Jacobi identity, and this 
implies that $\omega$ is a closed 2-form ($d\omega=0$). This $\omega$ 
is the symplectic form. Although right inverses are not unique, the 
Poisson bracket and symplectic form are independent of the specific 
choice of the $I^{\text{inv}}_N$'s here.

\subsection{Why coadjoint orbits?}
\label{coadjoint}
Let $\M$ be a compact manifold and assume that a compact, semisimple 
Lie group $G$ acts smoothly and transitively on $\M$. This implies 
that $\pi_1(\M)$ is finite, and thus $H^1(\M;\R)=0$. Everything we do 
should be $G$-equivariant. 

Because $G$ acts smoothly on $\M$, the elements of the Lie algebra 
$\g$ of $G$ define certain vector fields on $\M$. Since the 
quantization is assumed to be $G$-equivariant, the symplectic form 
must be $G$-invariant. This implies that for any $\xi\in\g$ thought 
of as a vector field on $\M$, 
$$
0=\Lie_\xi \omega = d(\xi\inner\omega)+\xi\inner d\omega = 
d(\xi\inner\omega)
\;;$$ 
so (using $H^1=0$) there is a ``Hamiltonian'' $h(\xi)\in\C^\infty(\M)$ such 
that $\xi\inner\omega=dh(\xi)$, which is well defined modulo 
constants. The constant can be fixed by requiring that the average of 
$h(\xi)$ over $\M$ is $0$. This gives a well-defined linear map 
$h:\g\to\C(\M)$.
For any $x\in\M$, the evaluation $\xi\mapsto h(\xi)(x)$ is a linear 
map $\g\to\co$; in other words, $h$ lets us map $x$ into the linear 
dual $\g^*$. That map is the ``moment map'' $\Phi:\M\to\g^*$ (see 
\cite{g-s}).

Because $\M$ is homogeneous and compact, the moment map turns out to 
be an embedding, so effectively $\M\subset\g^*$. By transitivity, 
$\M$ is precisely the orbit of any of its points under the natural 
``coadjoint'' action $\Ad^*_G$ of $G$ on $\g^*$, so any of the 
homogeneous spaces we are considering is a coadjoint orbit.

Since $\g^{**}=\g$, any element of $\g$ is naturally thought of as a 
linear function on $\g^*$. The Lie bracket on these of course 
satisfies the Jacobi identity, and extends to a unique Poisson 
bracket for all functions on $\g^*$. If $x_i$ are linear coordinates 
on $\g^*$ then the Poisson bivector $\pi$ on $\g^*$ is given by 
$\pi_{ij}=C^k_{\;ij}x_k$. This Poisson bivector is degenerate, but 
restricts to a nondegenerate one on any coadjoint orbit. This makes 
any coadjoint orbit symplectic. The set of homogeneous spaces we 
are interested in is therefore precisely the set of coadjoint orbits 
of compact Lie groups.

The single point $\{0\}\subset\g^*$ is trivially a coadjoint orbit. 
It is an exception to some of the statements in this paper, but an 
utterly uninteresting one, so I will not mention it again.

\subsection{Structure of coadjoint orbits}
\label{coadjoint_classification}
We are interested in all coadjoint orbits, but all coadjoint orbits 
intersect $\Ca^*\subset\g^*$, so it is sufficient to consider the 
orbits of all $\Lambda\in\Ca^*$. These are still not all distinct; 
$\Orl=\mathcal{O}_{\Lambda'}$ if (and only if) $\Lambda$ and 
$\Lambda'$ are mapped to one another by the Weyl group $W$. The set 
of distinct coadjoint orbits is therefore $\g^*/G \cong \Ca^*/W \cong 
\C_+$, using the fact that the Weyl chamber $\C_+$ is a fundamental 
domain of the $W$ action (App.~\ref{representation_theory}).

We would like to express the coadjoint orbit $\Orl$ of 
$\Lambda\in\C_+$ as $G/H$. So what is $H$? It is the subgroup of $G$ 
leaving $\Lambda$ invariant, or equivalently the centralizer 
\beq
H=Z_\Lambda (G) \equiv \{h\in G\mid \Ad_h(\Lambda)=\Lambda\}
\eeq
if $\Lambda$ is identified to an element of $\g$ using the 
inner product. In this sense, $\Lambda\in\Ca$, so because $\Ca$ 
is Abelian, $\Ca\subset\h$. This implies that the Cartan subgroup 
$T$ is a subgroup of $H$, so $T$ can be used as the Cartan subgroup 
of $H$, and weights of $G$ and $H$ are naturally identified. There are, 
however, weights which are dominant for $H$ that are not for $G$, and 
the Weyl groups are different.

Expand $\Lambda$ in the basis $\{\pi_j\}$, and mark the vertices 
$j\in\Dynkin(\g)$ for which $\pi_j$ has a nonzero coefficient in 
$\Lambda$. The vertices of the Dynkin diagram are also the natural 
index set for the dual basis of fundamental roots. In the standard 
root decomposition of $\g^\co$, $E_\alpha$ commutes with $\Lambda$ 
and is thus in $\h^\co$ if and only if $\alpha$ is orthogonal to 
$\Lambda$. This is true precisely if, in the expansion in fundamental 
roots, $\alpha$ has $0$ coefficients for all the marked vertices of 
$\Dynkin(\g)$. This means that $\h^\co$ is spanned by $\Ca$ and the 
$E_\alpha$'s that are supported on the unmarked vertices. 

This gives a simple, diagrammatic way of calculating $\h$: \emph{The 
Lie algebra $\h$ of $H$ is the sum of a copy of $\mathfrak{u}(1)$ for 
every marked vertex and the Lie algebra of whatever Dynkin diagram is 
left after deleting all the marked vertices (and adjoining edges).}

(This diagrammatic method is also described in \cite{b-f-r1}, the only 
difference is that the complementary set of vertices is marked.)

This shows that, up to homeomorphism, the orbit $\Orl$ depends only 
on which coefficients are nonzero. On the other hand the symplectic 
structure and metric do vary with $\Lambda$. Since the number of 
marked vertices is the number of nonzero coefficients for $\Lambda$, 
this is the number of parameters that orbits in a given homeomorphism 
class vary by. One of these degrees of freedom simply corresponds to 
rescaling.

I use the notation $[\lambda]$ for the irreducible $H$-representation 
with highest weight $\lambda$. If the weight $\lambda$ is a 
combination of fundamental weights corresponding to marked vertices 
of $\Dynkin(\g)$, then the semisimple part of $\h$ acts trivially on 
$[\lambda]$. In this case $[\lambda]$ is one-dimensional and is just 
a representation of the abelian part of $\h$. The weights $N\Lambda$ 
are of this type. 

In general, if $[\lambda]$ is one-dimensional then 
$[\lambda]^*=[-\lambda]$. If $[\lambda]$ is one-dimensional and $\mu$ 
is arbitrary, then $[\lambda]\otimes[\mu] = [\lambda+\mu]$.

\subsection{Examples}
\label{coadj_ex}
The existence of the symplectic structure implies that a coadjoint 
orbit must be even-dimensional. The table lists the 
lowest-dimensional coadjoint orbits (all those with dimension 
$\leq6$). Note that $\co P^3$ occurs in two forms.
\begin{table} 
\renewcommand{\arraystretch}{1.25}
\begin{center}
Coadjoint Orbits with $\dim\leq6$\\[1.25ex]
\begin{tabular}{|ccccc|}\hline
dim. & Name & $\Orl$ & $G/H$ & Diagram\\ \hline
2 & Sphere & $S^2$ & $\sut/\Uo$ & $\bullet$\\
4 & Complex projective space & $\co\mathrm{P}^2$ & 
$\mathrm{SU}(3)/\mathrm{U}(2)$ & 
$\mathrel\bullet\!\relbar\mkern-3.5mu\relbar\!\mathrel\circ$\\
6 & Complex projective space & $\co\mathrm{P}^3$ & 
$\mathrm{SU}(4)/\mathrm{U}(3)$ & 
$\mathrel\bullet\!\relbar\mkern-3.5mu\relbar\!\mathrel\circ\!\relbar
\mkern-3.5mu\relbar\!\mathrel\circ$ 
\\
& " & " &
$\mathrm{Sp}(4)/\Uo\times\mathrm{Sp}(2)$ &
$\mathrel\circ\!=\mkern-7.0mu\Rightarrow\!\mathrel\bullet$\\
6 & Complex flag variety & & $\mathrm{SU}(3)/\Uo\times\Uo$ & 
$\mathrel\bullet\!\relbar\mkern-3.5mu\relbar\!\mathrel\bullet$\\
6 & Double cover of real Grassmanian & $\widetilde G^\R_{2,5}$ & 
$\SO(5)/\SO(2)\times\SO(3)$ & $\mathrel\bullet 
\!=\mkern-7.0mu\Rightarrow\!\mathrel\circ$\\
 \hline
\end{tabular}
\end{center}\end{table}

A less trivial example is given by the diagram 
$\mathrel\bullet\!\relbar\mkern-3.5mu\relbar\!\mathrel\circ\!\relbar
\mkern-3.5mu\relbar\!\mathrel\bullet\!\relbar\mkern-3.5mu\relbar\!
\mathrel\circ\!=\mkern-7.0mu\Rightarrow\!\mathrel\circ$. 
In this case $G=\widetilde{\SO}(11)$, and (modulo coverings) 
$H\approx\Uo\times\sut\times\Uo\times\SO(5)$. The dimension is $\dim 
G - \dim H = 55-(1+3+1+10)=40$. 

Notably, $S^2$ is the only sphere which is a coadjoint orbit. In fact 
it is the only sphere which admits a symplectic structure, 
equivariant or not. This is because the symplectic form on a compact 
manifold always has a nontrivial cohomology class, implying 
$H^2(\M)\neq0$. The 2-sphere is the only sphere such that 
$H^2(S^n)\neq0$. This means that with the reasonable seeming 
condition of respecting the smooth structure (as described in 
\ref{symplectic}), no other sphere may be quantized. For a claim to 
the contrary, see \cite{g-k-p5}.

\section{Projective space}
\label{projective}
There is a (very) slightly different perspective on how the formula 
\eqref{I_map} for the injection $I_N:\A_N\into\C(\Orl)$ comes about. 
It can be thought of as resulting from a natural embedding of $\Orl$ 
into the projectivisation $\PP(N\Lambda)$ of the representation 
$(N\Lambda)$. The idea is simply that since $\Psi^{N\Lambda}$ is 
fixed modulo phase by $H$, its projective equivalence class 
$[\Psi^{N\Lambda}]\in \PP(N\Lambda)$ is exactly fixed by $H$. Indeed 
$H$ is the entire isotropy group of this point. This means that the 
equivariant map that takes $\Orl\ni o\mapsto[\Psi^{N\Lambda}]$ is an 
embedding $\Orl\into\PP(N\Lambda)$. 

Any point $[\psi]\in\PP(N\Lambda)$ determines a state (a normalized 
element of the dual) of $\A_N$. This takes  
\beq
a\mapsto \frac{\langle\psi|a|\psi\rangle}{\langle\psi|\psi\rangle} 
\;.\eeq
So, we can naturally map $\Orl\into\PP(N\Lambda)\to\A_N^*$. From this 
point the story continues in the same way as in Section 
\ref{coherent}.
\gap

Suppose that $\lvert x\rangle$, $\lvert y\rangle$, and $\lvert 
z\rangle$ are (unnormalized) vectors in $(\Lambda)$ such that 
$x\mapsto\left[\lvert x\rangle \right]\in\PP(\Lambda)$, \emph{et 
cetera}. The formula \eqref{kernel} for $K_1$ can be rewritten as 
\beq
\left[\dim(\Lambda)\right]^{-2} K_1(x,y,z) = 
\frac{\left<x|y\right>\left<y|z\right>\left<z|x\right>}
{\left<x|x\right>\left<y|y\right>\left<z|z\right>}
\label{K_unnorm}\;.\eeq
A continuous choice of these vectors cannot be made globally, but it 
can be made  in a small neighborhood of $o$. As in Section 
\ref{convergence}, let's fix $x=o$. The obvious choice for $\lvert 
o\rangle$ is $\lvert\Psi^\Lambda\rangle$. The arbitrariness in the 
other vectors is the freedom to multiply by a scalar. If we to fix 
these vectors by letting $\langle\Psi^\Lambda|y\rangle=1$ (and 
likewise for $z$), then \eqref{K_unnorm} simplifies to
\beq
\left[\dim (\Lambda)\right]^{-2} K_1(o,y,z) = 
\frac{\langle\Psi^\Lambda|y\rangle \langle y|z\rangle \langle 
z|\Psi^\Lambda\rangle}
{\langle\Psi^\Lambda|\Psi^\Lambda\rangle \langle y|y\rangle \langle 
z|z\rangle}
= \frac{\langle y|z\rangle}{\langle y|y\rangle \langle z|z\rangle}
\label{simp_K}\;.\eeq

Now suppose that we have a complex coordinate system for $y$ and $z$, 
that the coordinates are vectors $\upsilon$ and $\zeta$ in a subspace 
of $(\Lambda)$, and that to first order $\lvert y\rangle$ is given by 
$$
\lvert y\rangle \approx \lvert\Psi^\Lambda\rangle + 
\lvert\upsilon\rangle
$$
(and $\lvert z\rangle$ is given by $\zeta$). From this, the inner 
product $\langle y|z\rangle$ can be calculated to \emph{second} order
\begin{align*}
\langle y|z\rangle &= 1 + \langle y | z - \Psi^\Lambda\rangle
\\&
= 1 + \langle y-\Psi^\Lambda | z - \Psi^\Lambda\rangle
\\&
\approx 1 + \langle\upsilon|\zeta\rangle
\;.\end{align*}
Inserting this into \eqref{simp_K} gives a formula for $K_1$ to 
second order
\begin{align}
\left[\dim (\Lambda)\right]^{-2} K_1(o,y,z) 
&\approx \frac{1+ 
\langle\upsilon|\zeta\rangle}{\left(1+\norm{\upsilon}^2\right)\left(1+
\norm{\zeta}^2\right)}
\nonumber\\&
\approx 1 - \norm{\upsilon}^2 - \norm{\zeta}^2 + 
\langle\upsilon|\zeta\rangle
\label{K_approx}\;.\end{align}

\end{appendix}

\filbreak  
\section*{Notation}
\begin{tabular}{rl}
$[\:\cdot\:,\:\cdot\:]_-$ & Commutator, $[a,b]_-=ab-ba$. 
\eqref{J_commutation}.\\
$\A_{\hat\I}$ & The bundle of algebras over $\hat\I$. 
Sec.~\ref{generalQ}.\\
$\A_N$ & Quantum algebra at index $N\in\I$, later $\End(N\Lambda)$. 
Sec's \ref{generalQ}, \ref{Q_co_orb}.\\
$\AA$ & \cs-algebra of continuous sections of $\A_{\hat\I}$. Sec's 
\ref{generalQ}, \ref{gen_pic}, \ref{alg_dir_lim}, \ref{alg_inv_lim}.\\
$\oAA$ & Preliminary vector space, dense in $\AA$. \eqref{oAA_def}.\\
$\AA_0$ & $=\Gn(\A_\I)$, an ideal in $\AA$, the Kernel of 
$\Po:\AA\onto\Ai$. Sec's \ref{gen_pic}, \ref{alg_dir_lim}, 
\ref{alg_inv_lim}.\\ 
$\C_n$ & The $n$'th Casimir polynomial. \eqref{J_Casimir}.\\
$c_n(N\Lambda)$ & The eigenvalue of the Casimir operator $\C_n(J)$ 
acting on $(N\Lambda)$. \eqref{J_Casimir}.\\
$\Ca$ & Cartan subalgebra of $\g$. App.~\ref{representation_theory}.\\
$\C_+$ & Positive Weyl chamber in $\Ca^*$. Sec's 
\ref{intro_coadjoint}, \ref{representation_theory}.\\
$e_N$ & The function $e_N: \Orl\into\A_N,\; gH\mapsto 
\left|g\Psi^{N\Lambda}\right>\left<g\Psi^{N\Lambda}\right|$. 
\eqref{pformI}\\
$\End$ & Endomorphisms, the algebra of matrices on some vector space. 
Sec.~\ref{Q_co_orb}.\\
$\epsilon$ & Volume form on $\Orl$, normalized so that 
$\int_{\Orl}\epsilon=1$. \eqref{pformP}.\\
$\G$ & The space of continuous sections of a bundle. 
App.~\ref{sections}.\\
$\Gb$ & The space of norm-bounded sections. App.~\ref{sections}.\\
$\Gn$ & The space of continuous sections vanishing at $\infty$. 
App.~\ref{sections}.\\
$\Gc$ & The space of compactly supported sections. 
App.~\ref{sections}.\\
$\Gp$ & The space of polynomial sections. Sec.~\ref{twosphere}.\\
$\I$ & Index set of the quantization. Sec.~\ref{generalQ}.\\
$\hat\I$ & $=\I\cup\{\infty\}$. Sec.~\ref{generalQ}.\\
$J_i$ & The basis of Hermitian generators of $\g$. 
Sec.~\ref{gen_pic}.\\
$\N$ & $=\{1,2,\ldots\}$.\\
$\hat\N$ & $=\{1,2,\ldots,\infty\}$.\\
$(\lambda)$ & $G$-representation space with highest weight $\lambda$. 
Sec's \ref{Q_co_orb}, \ref{representation_theory}.\\
$[\lambda]$ & $H$-representation space with highest weight $\lambda$. 
Sec's \ref{identification}, \ref{coadjoint_classification}.\\
$\lambda^*$ & Weight vector such that $(\lambda^*)=(\lambda)^*$. 
Sec's \ref{coherent}, \ref{representation_theory}.\\
$\Orl$ & Coadjoint orbit passing through the weight vector 
$\Lambda\in\C_+\subset\g^*$. Sec.~\ref{intro_coadjoint}.\\
$\Po$ & The surjection $\AA\onto\C(\M)$. Sec.~\ref{generalQ}.\\
$\PP$ & Projectivization of a vector space. App.~\ref{projective}.\\
$\Pi$ & Projection onto some subrepresentation. Sec's \ref{coherent}, 
\ref{limit_bun}.\\
$\Psi^\lambda$ & Normalized highest weight vector in $(\lambda)$. 
Sec's \ref{coherent}, \ref{representation_theory}.\\
$Q^\lambda_N$ & Projection such that $V^\lambda_N = 
[\A_N\otimes(\mu)]Q_N$. Sec.~\ref{modules}\\
$\ntr$ & Trace normalized so that $\ntr ~1=1$. Sec.~\ref{coherent}.\\
$V_N$ & Module of the algebra $\A_N$ in the quantization of the 
bundle $V\!$. Sec.~\ref{vector}.\\
$V_\infty$ & $=\G(V)$, which is a module of $\C(\M)$. 
Sec.~\ref{vector}.\\
$V^\lambda_\N$ & The bundle of quantum modules associated to the 
$G$-weight $\lambda$. Sec.~\ref{modules}.\\
$\VV$ & $\AA$-module expressing a quantization of $V$\@. 
Sec.~\ref{vector}.\\
$\oVV$ & Vector space which generates $\VV$ as an $\AA$-module. App's 
\ref{mod_dir_lim}, \ref{mod_inv_lim}.\\
\end{tabular}

\subsection*{Acknowledgments}
I wish to thank Ranee Brylinski and Nigel Higson for extensive 
discussions. This material is based upon work supported under a 
National Science Foundation Graduate Fellowship. Also supported in 
part by NSF grant PHY95-14240 and by the Eberly Research Fund of the 
Pennsylvania State University

\end{document}